# HIGH-DIMENSIONAL ANALYSIS OF SEMIDEFINITE RELAXATIONS FOR SPARSE PRINCIPAL COMPONENTS[1]

By Arash A. Amini and Martin J. Wainwright

*University of California, Berkeley*

Principal component analysis (PCA) is a classical method for dimensionality reduction based on extracting the dominant eigenvectors of the sample covariance matrix. However, PCA is well known to behave poorly in the "large $p$, small $n$" setting, in which the problem dimension $p$ is comparable to or larger than the sample size $n$. This paper studies PCA in this high-dimensional regime, but under the additional assumption that the maximal eigenvector is sparse, say, with at most $k$ nonzero components. We consider a spiked covariance model in which a base matrix is perturbed by adding a $k$-sparse maximal eigenvector, and we analyze two computationally tractable methods for recovering the support set of this maximal eigenvector, as follows: (a) a simple diagonal thresholding method, which transitions from success to failure as a function of the rescaled sample size $\theta_{\mathrm{dia}}(n, p, k) = n/[k^2 \log(p-k)]$; and (b) a more sophisticated semidefinite programming (SDP) relaxation, which succeeds once the rescaled sample size $\theta_{\mathrm{sdp}}(n, p, k) = n/[k \log(p-k)]$ is larger than a critical threshold. In addition, we prove that no method, including the best method which has exponential-time complexity, can succeed in recovering the support if the order parameter $\theta_{\mathrm{sdp}}(n, p, k)$ is below a threshold. Our results thus highlight an interesting trade-off between computational and statistical efficiency in high-dimensional inference.

**1. Introduction.** Principal component analysis (PCA) is a classical method [1, 22] for reducing the dimension of data, say, from some high-dimensional subset of $\mathbb{R}^p$ down to some subset of $\mathbb{R}^d$, with $d \ll p$. Principal component analysis operates by projecting the data onto the $d$ directions of maximal

Received March 2008; revised August 2008.

[1]Supported in part by NSF Grants CAREER-CCF-05-45862 and DMS-06-05165, and a Sloan Foundation Fellowship.

variance, as captured by eigenvectors of the $p \times p$ population covariance matrix $\Sigma$. Of course, in practice, one does not have access to the population covariance, but instead must rely on a "noisy" version of the form

$$(1) \qquad \widehat{\Sigma} = \Sigma + \Delta,$$

where $\Delta = \Delta_n$ denotes a random noise matrix, typically arising from having only a finite number $n$ of samples. A natural question in assessing the performance of PCA is under what conditions the sample eigenvectors (i.e., based on $\widehat{\Sigma}$) are consistent estimators of their population analogues. In the classical theory of PCA, the model dimension $p$ is viewed as fixed, and asymptotic statements are established as the number of observations $n$ tends to infinity. With this scaling, the influence of the noise matrix $\Delta$ dies off, so that sample eigenvectors and eigenvalues are consistent estimators of their population analogues [1]. However, such "fixed $p$, large $n$" scaling may be inappropriate for many contemporary applications in science and engineering (e.g., financial time series, astronomical imaging, sensor networks), in which the model dimension $p$ is comparable or even larger than the number of observations $n$. This type of high-dimensional scaling causes dramatic breakdowns in standard PCA and related eigenvector methods, as shown by classical and ongoing work in random matrix theory [13, 20, 21].

Without further restrictions, there is little hope of performing high-dimensional inference with very limited data. However, many data sets exhibit additional structure, which can partially mitigate the curse of dimensionality. One natural structural assumption is that of sparsity, and various types of sparse models have been studied in past statistical work. There is a substantial and on-going line of work on subset selection and sparse regression models (e.g., [6, 11, 28, 35, 36]), focusing in particular on the behavior of various $\ell_1$-based relaxation methods. Other work has tackled the problem of estimating sparse covariance matrices in the high-dimensional setting, using thresholding methods [3, 12] as well as $\ell_1$-regularization methods [8, 39].

A related problem—and the primary focus of this paper—is recovering sparse eigenvectors from high-dimensional data. While related to sparse covariance estimation, the sparse eigenvector problem presents a different set of challenges; indeed, a covariance matrix may have a sparse eigenvector with neither it nor its inverse being a sparse matrix. Various researchers have proposed methods for extracting sparse eigenvectors, a problem often referred to as sparse principal component analysis (SPCA). Some of these methods are based on greedy or nonconvex optimization procedures (e.g., [23, 29, 40]), whereas others are based on various types of $\ell_1$-regularization [9, 41]. Zou, Hastie and Tibshirani [41] develop a method based on transforming the PCA problem to a regression problem and then applying the Lasso ($\ell_1$-regularization). Johnstone and Lu [21] proposed a two-step method, using an initial pre-processing step to select relevant variables followed by ordinary PCA in the reduced space. Under a particular $\ell_q$-ball sparsity model, they proved $\ell_2$-consistency of their procedure as long as $p/n$ converges to a constant. In recent work, d'Aspremont et al. [9] have formulated a direct semidefinite programming (SDP) relaxation of the sparse eigenvector problem, and developed fast algorithms for solving it, but have not provided high-dimensional consistency results. The elegant work of Paul and Johnstone [30, 32], brought to our attention after initial submission, studies estimation of eigenvectors satisfying weak $\ell_q$-ball sparsity assumptions for $q \in (0, 2)$. We discuss connections to this work at more length below.

In this paper, we study the model selection problem for sparse eigenvectors. More precisely, we consider a spiked covariance model [20], in which the maximal eigenvector $z^*$ of the population covariance $\Sigma_p \in \mathbb{R}^{p \times p}$ is $k$-sparse, meaning that it has nonzero entries on a subset $S(z^*)$ with cardinality $k$, and our goal is to recover this support set exactly. In order to do so, we have access to a matrix $\widehat{\Sigma}$, representing a noisy version of the population covariance, as in (1). Although our theory is somewhat more generally applicable, the most natural instantiation of $\widehat{\Sigma}$ is as a sample covariance matrix based on $n$ i.i.d. samples drawn from the population. We analyze this setup in the high-dimensional regime, in which all three parameters—the number of observations $n$, the ambient dimension $p$ and the sparsity index $k$—are allowed to tend to infinity simultaneously. Our primary interest is in the following question: using a given inference procedure, under what conditions on the scaling of triplet $(n, p, k)$ is it possible, or conversely impossible, to recover the support set of the maximal eigenvector $z^*$ with probability one?

We provide a detailed analysis of two procedures for recovering sparse eigenvectors, as follows: (a) a simple diagonal thresholding method, used as a pre-processing step by Johnstone and Lu [21], and (b) a semidefinite programming (SDP) relaxation for sparse PCA, recently developed by d'Aspremont et al. [9]. Under the $k$-sparsity assumption on the maximal eigenvector, we prove that the success or failure probabilities of these two methods have qualitatively different scaling in terms of the triplet $(n, p, k)$. For the diagonal thresholding method, we prove that its success or failure is governed by the rescaled sample size

$$(2) \qquad \theta_{\mathrm{dia}}(n, p, k) := \frac{n}{k^2 \log(p - k)},$$

meaning that it succeeds with probability one for scalings of the triplet $(n, p, k)$ such that $\theta_{\mathrm{dia}}$ is above some critical value and, conversely, fails with probability one when this ratio falls below some critical value. We then establish performance guarantees for the SDP relaxation [9]. In particular, for the same class of models, we show that it always has a unique rank-one solution that specifies the correct signed support once $\theta_{\mathrm{dia}}(n, p, k)$ is sufficiently large, moreover, that for sufficiently large values of the rescaled sample size

$$(3) \qquad \theta_{\mathrm{sdp}}(n, p, k) := \frac{n}{k \log(p - k)},$$

if there exists a rank-one solution, then it specifies the correct signed support. The proof of this result is based on random matrix theory, concentration of measure and Gaussian comparison inequalities. Our final contribution is to use information-theoretic arguments to show that no method can succeed in recovering the signed support for the spiked identity covariance model if the order parameter $\theta_{\mathrm{sdp}}(n, p, k)$ lies below some critical value. One consequence is that the given scaling (3) for the SDP relaxation is sharp, meaning the SDP relaxation also fails once $\theta_{\mathrm{sdp}}$ drops below a critical threshold. Moreover, it shows that under the rank-one condition, the SDP is in fact statistically optimal, that is, it requires only the necessary number of samples (up to a constant factor) to succeed.

The results reported here are complementary to those of Paul and Johnstone [30, 32], who propose and analyze the augmented SPCA algorithm for estimating eigenvectors. In comparison to the models analyzed here, their analysis applies to spiked models using the identity base covariance, but it allows for $m > 1$ eigenvectors in the spiking. In addition, they consider the class of weak $\ell_q$-ball sparsity models, as opposed to the hard $\ell_0$-sparsity model considered here. Another difference is that their results provide guarantees in terms of the $\ell_2$-norm between the eigenvector and its estimate, whereas our results guarantee exact support recovery. We note that an estimate can be close in $\ell_2$-norm while having a very different support set. Consequently, the results given here, which provide conditions for exact support recovery, provide complementary insight.

Our results highlight some interesting trade-offs between computational and statistical costs in high-dimensional inference. On one hand, the statistical efficiency of SDP relaxation is substantially greater than the diagonal thresholding method, requiring $\mathcal{O}(1/k)$ fewer observations to succeed. However, the computational complexity of SDP is also larger by roughly a factor $\mathcal{O}(p^3)$. An implementation due to d'Asprémont et al. [9] has complexity $\mathcal{O}(np + p^4 \log p)$ as opposed to the $\mathcal{O}(np + p \log p)$ complexity of the diagonal thresholding method. Moreover, our information-theoretic analysis shows that the best possible method—namely, one based on an exhaustive search over all $\binom{p}{k}$ subsets, with exponential complexity—does not have substantially greater statistical efficiency than the SDP relaxation.

The remainder of this paper is organized as follows. In Section 2, we provide precise statements of our main results, discuss some of their implications and provide simulation results to illustrate the sharpness of their predictions. Sections 3, 4 and 5 are devoted to proofs of these results, with some of the more technical aspects deferred to appendices. We conclude in Section 6.

1.1. *Notation.* For the reader's convenience, we state here some notation used throughout the paper. For a vector $x \in \mathbb{R}^n$, we use $\|x\|_p = (\sum_{i=1}^n |x_i|^p)^{1/p}$ to denote its $\ell_p$-norm. For a matrix $A \in \mathbb{R}^{m \times n}$, we use $\|A\|_{p,q}$ to denote the matrix operator norm induced by vector norms $\ell_p$ and $\ell_q$; more precisely, we have

$$(4) \qquad \|A\|_{p,q} := \max_{\|x\|_q = 1} \|Ax\|_p.$$

A few cases of particular interest in this paper are (a) the *spectral norm* given by

$$\|A\|_{2,2} := \max_{i=1,\ldots,m} \{\sigma_i(A)\},$$

where $\{\sigma_i(A)\}$ are the singular values of $A$, and the $\ell_\infty$-*operator norm*, given by

$$\|A\|_{\infty,\infty} := \max_{i=1,\ldots,m} \sum_{j=1}^n |A_{ij}|.$$

Given two square matrices $X, Y \in \mathbb{R}^{n \times n}$, we define the matrix inner product $\langle\!\langle X, Y \rangle\!\rangle := \operatorname{tr}(XY^T) = \sum_{i,j} X_{ij} Y_{ij}$. Note that this inner product induces the Hilbert–Schmidt norm $\|X\|_{\mathrm{HS}} = \sqrt{\langle\!\langle X, X \rangle\!\rangle}$.

We use the following standard asymptotic notation: for functions $f, g$, the notation $f(n) = \mathcal{O}(g(n))$ means that there exists a fixed constant $0 < C < +\infty$ such that $f(n) \leq C g(n)$; the notation $f(n) = \Omega(g(n))$ means that $f(n) \geq C g(n)$, and $f(n) = \Theta(g(n))$ means that $f(n) = \mathcal{O}(g(n))$ and $f(n) = \Omega(g(n))$. Note in particular that when used without a subscript "$p$," these symbols are to be interpreted in a deterministic sense, that is, the constants involved are assumed to be nonrandom.

We use $\lambda(A)$ to denote a generic eigenvalue of a square matrix $A$, as well as $\lambda_{\min}(\cdot)$ and $\lambda_{\max}(\cdot)$ for the minimal and the maximal eigenvalues, respectively. Any member of the set of eigenvectors of $A$ associated with an eigenvalue is denoted as $\vec{v}(A)$. Thus, $\vec{v}_{\max}(\cdot)$, for example, represents the eigenvectors associated with the maximal eigenvalue (occasionally referred to as "maximal eigenvectors"). We always assume that eigenvectors are normalized to unit $\ell_2$-norm and have a nonnegative first component. The sign convention guarantees uniqueness of the eigenvector associated with an eigenvalue with geometric multiplicity one.

Finally, some probabilistic notation: we say a sequence of events $\{E_j\}_{j \geq 1}$ happens with asymptotic probability one (w.a.p. one) if $\lim_{j \to +\infty} \mathbb{P}[E_j] = 1$, whereas it holds asymptotically almost surely (a.a.s.) as $j \to +\infty$ if $\mathbb{P}(\liminf E_j) = 1$.

**2. Main results and consequences.** The primary focus of this paper is the *spiked covariance model*, in which some base covariance matrix is perturbed by the addition of a sparse eigenvector $z^* \in \mathbb{R}^p$. In particular, we study sequences of covariance matrices of the form

$$(5) \qquad \Sigma_p = \beta z^* z^{*T} + \begin{bmatrix} I_k & 0 \\ 0 & \Gamma_{p-k} \end{bmatrix} = \beta z^* z^{*T} + \Gamma,$$

where $\Gamma_{p-k} \in \mathbb{S}_+^{p-k}$ is a symmetric PSD matrix with $\lambda_{\max}(\Gamma_{p-k}) \le 1$. Note that we have assumed (without loss of generality, by re-ordering the indices as necessary) that the nonzero entries of $z^*$ are indexed by $\{1, \ldots, k\}$, so that (5) is the form of the covariance after any re-ordering. We also assume that the nonzero part of $z^*$ has entries $z_i^* \in \frac{1}{\sqrt{k}}\{-1, +1\}$, so that $\|z^*\|_2 = 1$.

The spiked covariance model (5) was first proposed by Johnstone [20], who focused on the spiked identity covariance matrix [i.e., model (5) with $\Gamma_{p-k} = I_{p-k}$]. Johnstone and Lu [21] established that the sample eigenvectors for the spiked identity model, based on a set of $n$ i.i.d. samples with distribution $N(0, \Sigma_p)$ from the spiked identity ensemble, are inconsistent as estimators of $z_i^*$ whenever $p/n \to c > 0$. These asymptotic results were refined by later work [2, 31].

In this paper, we study a slightly more general family of spiked covariance models, in which the matrix $\Gamma_{p-k}$ is required to satisfy the following conditions:

$$(6a) \qquad \text{A1.} \qquad \|\sqrt{\Gamma_{p-k}}\|_{\infty,\infty} = \mathcal{O}(1) \quad \text{and}$$

$$(6b) \qquad \text{A2.} \qquad \lambda_{\max}(\Gamma_{p-k}) \le \min\left\{1, \lambda_{\min}(\Gamma_{p-k}) + \frac{\beta}{8}\right\}.$$

Here $\sqrt{\Gamma_{p-k}}$ denotes the symmetric square root. These conditions are trivially satisfied by the identity matrix $I_{p-k}$, but also can hold for more general nondiagonal matrices. Thus, under the model (5), the population covariance matrix $\Sigma$ itself need not be sparse, since (at least generically) it has $k^2 + (p-k)^2 = \Theta(p^2)$ nonzero entries. Assumption (A2) on the eigenspectrum of the matrix $\Gamma_{p-k}$ ensures that as long as $\beta > 0$, then the vector $z^*$ is the unique maximal eigenvector of $\Sigma$, with associated eigenvalue $(1 + \beta)$. Since the remaining eigenvalues are bounded above by 1, the parameter $\beta > 0$ represents a signal-to-noise ratio, characterizing the separation between the maximal eigenvalue and the remainder of the eigenspectrum. Assumption (A1) is related to the fact that recovering the correct signed support means that the estimate $\widehat{z}$ must satisfy $\|\widehat{z} - z^*\|_\infty \le 1/\sqrt{k}$. As will be clarified by our analysis (see Section 4.4), controlling this $\ell_\infty$-norm requires bounds on terms of the form $\|\sqrt{\Gamma_{p-k}}u\|_\infty$, which requires control of the $\ell_\infty$-operator norm $\|\sqrt{\Gamma_{p-k}}\|_{\infty,\infty}$.

In this paper, we study the *model selection problem* for eigenvectors: that is, we assume that the maximal eigenvector $z^*$ is $k$-sparse, meaning that it has exactly $k$ nonzero entries, and our goal is to recover this support, along with the sign of $z^*$ on its support. We let $S(z^*) = \{i \mid z_i^* \neq 0\}$ denote the support set of the maximal eigenvector; recall that $S(z^*) = \{1, \ldots, k\}$ by our assumed ordering of the indices. Moreover, we define the function $S_\pm : \mathbb{R}^p \to \{-1, 0, +1\}^p$ by

$$[S_\pm(u)]_i := \begin{cases} \text{sign}(u_i), & \text{if } u_i \neq 0, \\ 0, & \text{otherwise,} \end{cases} \tag{7}$$

so that $S_\pm(z^*)$ encodes the *signed support* of the maximal eigenvector.

Given some estimate $\widehat{S_\pm}$ of the true signed support $S_\pm(z^*)$, we assess it based on the 0–1 loss $\mathbb{I}[\widehat{S_\pm} \neq S_\pm(z^*)]$, so that the associated risk is simply the probability of incorrect decision $\mathbb{P}[\widehat{S_\pm} \neq S_\pm(z^*)]$. Our goal is to specify conditions on the scaling of the triplet $(n, p, k)$ such that this error probability vanishes, or conversely, fails to vanish asymptotically. We consider methods that operate based on a set of $n$ samples $x^1, \ldots, x^n$, drawn i.i.d. with distribution $N(0, \Sigma_p)$. Under the spiked covariance model (5), each sample can be written as

$$x^i = \sqrt{\beta} v^i z^* + \sqrt{\Gamma} g^i, \tag{8}$$

where $\sqrt{\Gamma}$ is the symmetric matrix square root. Here $v^i \sim N(0, 1)$ is standard Gaussian, and $g^i \sim N(0, I_{p \times p})$ is a standard Gaussian $p$-vector, independent of $v^i$, so that $\sqrt{\Gamma} g^i \sim N(0, \Gamma)$. The data $\{x^i\}_{i=1}^n$ defines the sample covariance matrix

$$\widehat{\Sigma} := \frac{1}{n} \sum_{i=1}^n (x^i)(x^i)^T, \tag{9}$$

which follows a $p$-variate Wishart distribution [1]. In this paper, we analyze the high-dimensional scaling of two methods for recovering the signed support of the maximal eigenvector. It will be assumed throughout that the size $k$ of the support of $z^*$ is available to the methods a priori, that is, we do not make any attempt at estimating $k$.

2.1. *Diagonal thresholding method.* Under the spiked covariance model (5), note that the diagonal elements of the population covariance satisfy $\Sigma_{\ell\ell} = 1 + \beta/k$ for all $\ell \in S$, and $\Sigma_{\ell\ell} \leq 1$ for all $\ell \notin S$. (This latter bound follows since for all $\ell \notin S$, we have $\Sigma_{\ell\ell} \leq \|\Gamma_{p-k}\|_{2,2} \leq 1$.) This observation motivates a natural approach to recovering information about the support set $S$, previously used as a pre-processing step by Johnstone and Lu [21].

Let $D_\ell, \ell = 1, \ldots, p$, be the diagonal elements of the sample covariance matrix—namely,

$$D_\ell = \frac{1}{n} \sum_{i=1}^{n} (x_\ell^i)^2 = [\widehat{\Sigma}]_{\ell\ell}.$$

Form the associated order statistics

$$D_{(1)} \leq D_{(2)} \leq \cdots \leq D_{(p-1)} \leq D_{(p)},$$

and output the random subset $\widehat{S}(D)$ of cardinality $k$ specified by the indices of the largest $k$ elements $\{D_{(p-k+1)}, \ldots, D_{(p)}\}$. The chief appeal of this method is its low computational complexity. Apart from the order $\mathcal{O}(np)$ of computing the diagonal elements of $\widehat{\Sigma}$, it requires only performing a sorting operation, with complexity $\mathcal{O}(p \log p)$.

Note that this method provides only an estimate of the support $S(z^*)$, as opposed to the signed support $S_\pm(z^*)$. One could imagine extending the method to extract sign information as well, but our main interest in studying this method is to provide a simple benchmark by which to calibrate our later results on the performance of the more complex SDP relaxation. In particular, the following result provides a precise characterization of the statistical behavior of the diagonal thresholding method.

PROPOSITION 1 (Performance of diagonal thresholding). *For $k = \mathcal{O}(p^{1-\delta})$ for any $\delta \in (0, 1)$, the probability of successful recovery using diagonal thresholding undergoes a phase transition as a function of the rescaled sample size*

$$\theta_{\mathrm{dia}}(n, p, k) = \frac{n}{k^2 \log(p - k)}. \tag{10}$$

*More precisely, there exists a constant $\theta_u$ such that if $n > \theta_u k^2 \log(p - k)$, then*

$$\mathbb{P}[\widehat{S}(D) = S(z^*)] \geq 1 - \exp(-\Theta(k^2 \log(p - k))) \to 1, \tag{11}$$

*so that the method succeeds w.a.p. one and a constant $\theta_\ell > 0$ such that if $n \leq \theta_\ell k^2 \log(p - k)$, then*

$$\mathbb{P}[\widehat{S}(D) = S(z^*)] \leq \exp(-\Theta(\log(p - k))) \to 0, \tag{12}$$

*so that the method fails w.a.p. one.*

REMARKS. The proof of Proposition 1, provided in Section 3, is based on large deviations bounds on $\chi^2$-variates. The achievability assertion (11) uses known upper bounds on the tails of $\chi^2$-variates (e.g., [4, 21]). The converse result (12) requires an exponentially tight lower bound on the tails of $\chi^2$-variates, which we derive in Appendix C.

To illustrate the prediction of Proposition 1, we provide some results on the diagonal thresholding method. For all experiments reported here, we generated $n$ samples $\{x^1, \ldots, x^n\}$ in an i.i.d. manner from the spiked covariance ensemble (5), with $\Gamma = I$ and $\beta = 3$. Figure 1 illustrates the behavior predicted by Proposition 1. Each panel plots the success probability $\mathbb{P}[\widehat{S}(D) = S(z^*)]$ versus the rescaled sample size $\theta_{\mathrm{dia}}(n, p, n) = n/[k^2 \log(p - k)]$. Each panel shows five model dimensions ($p \in \{100, 200, 300, 600, 1200\}$), with panel (a) showing the logarithmic sparsity index $k = \mathcal{O}(\log p)$ and panel (b) showing the case $k = \mathcal{O}(\sqrt{p})$. Each point on each curve corresponds to the average of 100 independent trials. As predicted by Proposition 1, the curves all coincide, even though they correspond to very different regimes of $(p, k)$.

2.2. *Semidefinite-programming relaxation.* We now describe the approach to sparse PCA developed by d'Asprémont et al. [9]. Let $\mathbb{S}^p_+ = \{Z \in \mathbb{R}^{p \times p} \mid Z = Z^T, Z \succeq 0\}$ denote the cone of symmetric, positive semidefinite (PSD) matrices. Given $n$ i.i.d. observations from the model $N(0, \Sigma_p)$, let $\widehat{\Sigma}$ be the sample covariance matrix (9), and let $\rho_n > 0$ be a user-defined regularization parameter. d'Asprémont et al. [9] propose estimating $z^*$ by solving the optimization problem

$$(13) \qquad \widehat{Z} := \underset{Z \in \mathbb{S}^p_+}{\arg\max} \left[ \operatorname{tr}(\widehat{\Sigma} Z) - \rho_n \sum_{i,j} |Z_{ij}| \right] \qquad \text{s.t. } \operatorname{tr}(Z) = 1,$$

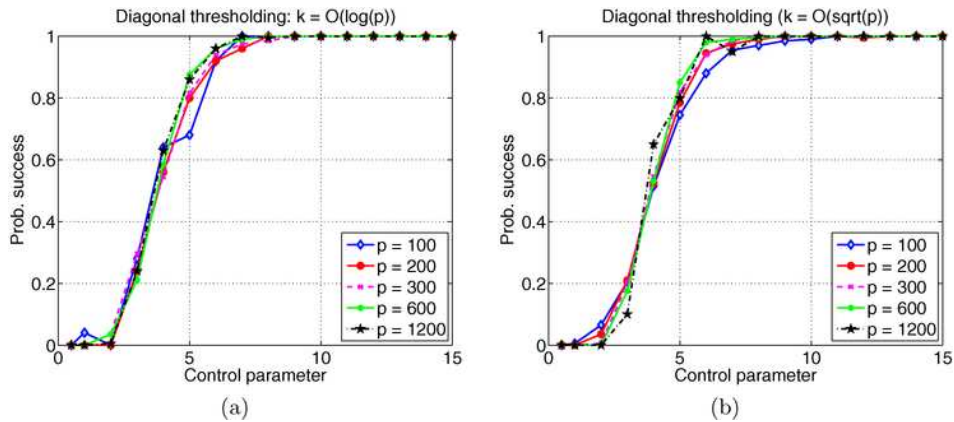

(a)             (b)

FIG. 1. *Plot of the success probability $\mathbb{P}[\widehat{S}(D) = S(z^*)]$ versus the rescaled sample size $\theta_{\mathrm{dia}}(n, p, k) = n/[k^2 \log(p - k)]$. The five curves in each panel correspond to model dimensions $p \in \{100, 200, 300, 600, 1200\}$, SNR parameter $\beta = 3$ and sparsity indices $k = \mathcal{O}(\log p)$ in panel (a) and $k = \mathcal{O}(\sqrt{p})$ in panel (b). As predicted by Proposition 1, the success probability undergoes a phase transition, with the curves for different model sizes and different sparsity indices all lying on top of one another.*

and computing the maximal eigenvector $\hat{z} = \vec{v}_{\max}(\widehat{Z})$. The optimization problem (13) is a semidefinite program (SDP), a class of convex conic programs that can be solved exactly in polynomial time. Indeed, d'Asprémont et al. [9] describe an $\mathcal{O}(p^4 \log p)$ algorithm, with an implementation posted online, that we use for all simulations reported in this paper.

To gain some intuition for the SDP relaxation (13), recall the following Courant–Fischer variational representation [18] of the maximal eigenvalue and eigenvector:

$$(14) \qquad \vec{v}_{\max}(\widehat{\Sigma}) = \arg\max_{\|z\|_2 = 1} z^T \widehat{\Sigma} z.$$

A lesser known but equivalent variational representation is in terms of the semidefinite program (SDP)

$$(15) \qquad Z^* = \arg\max_{Z \in \mathbb{S}^p_+, \mathrm{tr}(Z) = 1} \mathrm{tr}(\widehat{\Sigma} Z).$$

For this problem, if the maximal eigenvalue is simple, the optimum is always achieved at a rank-one matrix $Z^* = z^*(z^*)^T$, where $z^* = \vec{v}_{\max}(\widehat{\Sigma})$ is the maximal eigenvector; otherwise, there exist optimal solutions of higher rank, but the optimum is always achieved by at least some rank-one matrix. If we were given a priori information that the maximal eigenvector were sparse, then it might be natural to solve the same semidefinite program with the addition of an $\ell_0$ constraint. Given the intractability of such an $\ell_0$-optimization problem, the SDP program (13) is a natural relaxation.

In particular, the following result provides sufficient conditions for the SDP relaxation (13) to succeed in recovering the correct signed support of the maximal eigenvector.

THEOREM 2 (SDP performance guarantees). *Impose conditions (6a) and (6b) on the sequence of population covariance matrices $\{\Sigma_p\}$, and suppose moreover that $\rho_n = \beta/(2k)$ and $k = \mathcal{O}(\log p)$. Then:*

(a) Rank guarantee: *there exists a constant $\theta_{\mathrm{wr}} = \theta_{\mathrm{wr}}(\Gamma, \beta)$ such that for all sequences $(n, p, k)$ satisfying $\theta_{\mathrm{dia}}(n, p, k) > \theta_{\mathrm{wr}}$, the semidefinite program (13) has a rank-one solution with high probability.*

(b) Critical scaling: *there exists a constant $\theta_{\mathrm{crit}} = \theta_{\mathrm{crit}}(\Gamma, \beta)$ such that if the sequence $(n, p, k)$ satisfies*

$$(16) \qquad \theta_{\mathrm{sdp}}(n, p, k) := \frac{n}{k \log(p - k)} > \theta_{\mathrm{crit}}$$

*and if there exists a rank-one solution, then it specifies the correct signed support with probability converging to one.*

REMARKS. Part (a) of the theorem shows that rank-one solutions of the SDP (13) are not uncommon; in particular, they are guaranteed to exist with high probability at least under the weaker scaling of the diagonal thresholding method. The main contribution of Theorem 2 is its part (b), which provides sufficient conditions for signed support recovery using the SDP, when a rank-one solution exists. The bulk of our technical effort is devoted to part (b); indeed, the proof of part (a) is straightforward once all the pieces of the proof of part (b) have been introduced, and so will be deferred to Appendix G. For technical reasons, our current proof(s) require the condition $k = \mathcal{O}(\log p)$; however, it should be possible to remove this restriction, and indeed, the empirical results do not appear to require it.

Proposition 1 and Theorem 2 apply to the performance of specific (polynomial-time) methods. It is natural then to ask whether there exists any algorithm, possibly with super-polynomial complexity, that has greater statistical efficiency. The following result is information-theoretic in nature, and characterizes the fundamental limitations of any algorithm regardless of its computational complexity.

THEOREM 3 (Information-theoretic limitations). *Consider the problem of recovering the eigenvector support in the spiked covariance model (5) with $\Gamma = I_p$. For any sequence $(n, p, k) \to +\infty$ such that*

$$(17) \qquad \theta_{\mathrm{sdp}}(n, p, k) := \frac{n}{k \log(p - k)} < \frac{1 + \beta}{\beta^2},$$

*the probability of error of any method is at least $1/2$.*

REMARKS. Together with Theorem 2, this result establishes the sharpness of the threshold (16) in characterizing the behavior of SDP relaxation, and moreover, it guarantees optimality of the SDP scaling (16), up to constant factors, for the spiked identity ensemble.

To illustrate the predictions of Theorem 2 and 3, we applied the SDP relaxation to the spiked identity covariance ensemble, again generating $n$ i.i.d. samples. We solved the SDP relaxation using publically available code provided by d'Asprémont et al. [9]. Figure 2 shows the corresponding plots for the SDP relaxation [9]. Here we plot the probability $\mathbb{P}[S_\pm(\hat{z}) = S_\pm(z^*)]$ that the SDP relaxation correctly recovers the signed support of the unknown eigenvector $z^*$, where the signs are chosen uniformly in $\{-1, +1\}$ at random. Following Theorem 2, the horizontal axis plots the rescaled sample size $\theta_{\mathrm{sdp}}(n, p, k) = n/[k \log(p - k)]$. Each panel shows plots for three different problem sizes, $p \in \{100, 200, 300\}$, with panel (a) corresponding to logarithmic sparsity $[k = \mathcal{O}(\log p)]$, and panel (b) to linear sparsity $(k = 0.1p)$.

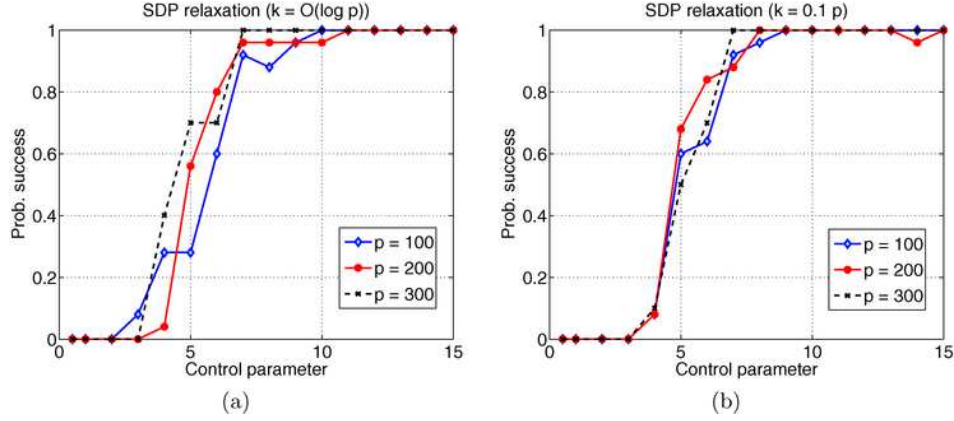

FIG. 2. *Performance of the SDP relaxation for the spiked identity ensemble, plotting the success probability $\mathbb{P}[S_\pm(\widehat{z}) = S_\pm(z^*)]$ versus the rescaled sample size $\theta_{\mathrm{sdp}}(n, p, k) = n/[k \log(p - k)]$. The three curves in each panel correspond to model dimensions $p \in \{100, 200, 300\}$, SNR parameter $\beta = 3$ and sparsity indices $k = \mathcal{O}(\log p)$ in panel* (a) *and $k = 0.1p$ in panel* (b). *As predicted by Theorem 2, the curves in panel* (a) *all lie on top of one another, and transition to success once the order parameter $\theta_{\mathrm{sdp}}$ is sufficiently large.*

Consistent with the prediction of Theorem 2, the success probability rapidly approaches one once the rescaled sample size exceeds some critical threshold. [Strictly speaking, Theorem 2 only covers the case of logarithmic sparsity shown in panel (a), but the linear sparsity curves in panel (b) show the same qualitative behavior.] Note that this empirical behavior is consistent with our conclusion that the order parameter $\theta_{\mathrm{sdp}}(n, p, k) = n/[k \log(p - k)]$ is a sharp description of the SDP threshold.

**3. Proof of Proposition 1.** We begin by proving the achievability result (11). We provide a detailed proof for the case $\Gamma_{p-k} = I_{p-k}$ and discuss necessary modifications for the general case at the end. For $\ell = 1, \dots, p$, we have

$$D_\ell = \frac{1}{n} \sum_{i=1}^{n} (x_\ell^i)^2 = \frac{1}{n} \sum_{i=1}^{n} [\sqrt{\beta} z_\ell^* v^i + g_\ell^i]^2. \tag{18}$$

Since $(\sqrt{\beta} z_\ell^* v^i + g_\ell^i) \sim N(0, \beta(z_\ell^*)^2 + 1)$ for each $i$, the rescaled variate $\frac{n}{\beta(z_\ell^*)^2 + 1} \times D_\ell$ is central $\chi_n^2$ with $n$ degrees of freedom. Consequently, we have

$$\mathbb{E}[D_\ell] = \begin{cases} 1, & \text{for all } \ell \in S^c, \\ 1 + \dfrac{\beta}{k}, & \text{for all } \ell \in S, \end{cases}$$

where we have used the fact that $(z_\ell^*)^2 = 1/k$ for $\ell \in S$, by assumption.

A sufficient condition for success of the diagonal thresholding decoder is a threshold $\tau_k$ such that $D_\ell \geq (1 + \tau_k)$ for all $\ell \in S$, and $D_\ell < (1 + \tau_k)$ for all $\ell \in S^c$. Using the union bound and the tail bound (61) on central $\chi^2$, we have

$$\mathbb{P}\Big[\max_{\ell \in S^c} D_\ell \geq (1 + \tau_k)\Big] \leq (p - k)\mathbb{P}\Big[\frac{\chi_n^2}{n} \geq 1 + \tau_k\Big] \leq (p - k)\exp\Big(-\frac{3n}{16}\tau_k^2\Big),$$

so that the probability of false inclusion vanishes as long as $n > \frac{16}{3}(\tau_k)^{-2} \times \log(p - k)$.

On the other hand, using the union bound and the tail bound (60b), we have

$$\mathbb{P}\Big[\min_{\ell \in S} D_\ell < (1 + \tau_k)\Big] \leq k\mathbb{P}\Big[\frac{\chi_n^2}{n} - 1 < \frac{1 + \tau_k}{1 + \beta/k} - 1\Big]$$

$$= k\mathbb{P}\Big[\frac{\chi_n^2}{n} - 1 < \frac{\tau_k - \beta/k}{1 + \beta/k}\Big]$$

$$\leq k\mathbb{P}\Big[\frac{\chi_n^2}{n} - 1 < \tau_k - \frac{\beta}{k}\Big].$$

As long as $\tau_k < \beta/k$, we may choose $x = \frac{n}{4}(\frac{\beta}{k} - \tau_k)^2$ in (60b), thereby obtaining the upper bound

$$\mathbb{P}\Big[\min_{\ell \in S} D_\ell < n(1 + \tau_k)\Big] \leq k\exp\Big(-\frac{n}{4}\Big(\frac{\beta}{k} - \tau_k\Big)^2\Big),$$

so that the probability of false exclusion vanishes as long as $n > \frac{4}{(\beta/k - \tau_k)^2}\log k$. Overall, choosing $\tau_k = \frac{\beta}{2k}$ ensures that the probability of both types of error vanish asymptotically as long as

$$n > \max\Big\{\frac{64}{3\beta^2}k^2\log(p - k), \frac{16}{\beta^2}k^2\log k\Big\}.$$

Since $k = o(p)$, the $\log(p - k)$ term is the dominant requirement. The modifications required for the case of general $\Gamma_{p-k}$ are straightforward. Since $\mathrm{var}(\sqrt{\Gamma}g^i)_\ell = (\Gamma_{p-k})_{\ell\ell} \leq 1$ for all $\ell \in S^c$ and samples $i = 1, \ldots, n$, we need to adjust the scaling of the $\chi_n^2$ variates. For general $\Gamma_{p-k}$, the variates $\{D_\ell, \ell \in S^c\}$ need no longer be independent, but our proof used only union bound, and so is valid regardless of the dependence structure.

We now prove the converse claim (12) for the spiked identity ensemble. At a high level, this portion of the proof consists of the following steps. For a positive real $t$, define the events

$$\mathbb{A}_1(t) := \Big\{\max_{\ell \in S^c} D_\ell > 1 + t\Big\} \quad \text{and} \quad \mathbb{A}_2(t) := \Big\{\min_{\ell \in S} D_\ell < 1 + t\Big\}.$$

Noting that the event $\mathbb{A}_1(t) \cap \mathbb{A}_2(t)$ implies failure of the diagonal cutoff decoder, it suffices to show the existence of some $t > 0$ such that $\mathbb{P}[\mathbb{A}_1(t)] \to 1$ and $\mathbb{P}[\mathbb{A}_2(t)] \to 1$.

*Analysis of event* $\mathbb{A}_1(t)$. Central to the analysis of event $\mathbb{A}_1$ is the following large-deviations *lower bound* on $\chi^2$-variates.

LEMMA 4. *For a central $\chi_n^2$ variable with $n$ degrees of freedom, there exists a constant $C > 0$ such that*

$$\mathbb{P}\Big[\frac{\chi_n^2}{n} > 1 + t\Big] \geq \frac{C}{\sqrt{n}} \exp(-nt^2/2)$$

*for all $t \in (0, 1)$.*

See Appendix C for the proof.

We exploit this lemma as follows. First, define the integer-valued random variable

$$Z(t) := \sum_{\ell \in S^c} \mathbb{I}[D_\ell > 1 + t]$$

corresponding to the number of indices $\ell \in S^c$ for which the diagonal entry $D_\ell$ exceeds $1 + t$, and note that $\mathbb{P}[\mathbb{A}_1(t)] = \mathbb{P}[Z(t) > 0]$. By a one-sided Chebyshev inequality [15], we have

$$(19) \qquad \mathbb{P}[\mathbb{A}_1(t)] = \mathbb{P}[Z(t) > 0] \geq \frac{(\mathbb{E}[Z(t)])^2}{(\mathbb{E}[Z(t)])^2 + \operatorname{var}(Z(t))}.$$

Note that $Z(t)$ is a sum of $(p - k)$ independent Bernoulli indicators, each with the same parameter $q(t) := \mathbb{P}[D_\ell > 1 + t]$. Computing the mean $\mathbb{E}[Z(t)] = (p - k)q(t)$ and variance $\operatorname{var}(Z(t)) = (p - k)q(t)(1 - q(t))$, and then substituting into the Chebyshev bound (19), we obtain

$$\mathbb{P}[\mathbb{A}_1(t)] \geq \frac{(p-k)^2 q^2(t)}{(p-k)^2 q^2(t) + (p-k)q(t)(1-q(t))} \geq \frac{(p-k)q(t)}{(p-k)q(t) + 1}$$

$$\geq 1 - \frac{1}{(p-k)q(t)}.$$

Consequently, the condition $(p - k)q(t) \to \infty$ implies that $\mathbb{P}[\mathbb{A}_1(t)] \to 1$.

Let us set $t = \sqrt{\frac{\delta \log(p-k)}{n}}$. [Here $\delta \in (0, 1)$ is the parameter from the assumption $k = \mathcal{O}(p^{1-\delta})$.] From Lemma 4, we have $q(t) \geq \frac{C}{\sqrt{n}} \exp(-nt^2/2)$, so that

$$(p-k)q\left(\sqrt{\frac{\delta \log(p-k)}{n}}\right) \geq \frac{C(p-k)}{\sqrt{n}} \exp\left(-\frac{\delta}{2} \log(p-k)\right)$$

$$= \frac{C(p-k)^{1-\delta/2}}{\sqrt{n}}.$$

Since $n \leq Lk^2 \log(p - k)$ for some $L < +\infty$ by assumption, we have

$$(p - k)q\left(\sqrt{\frac{\delta \log(p - k)}{n}}\right) \geq \frac{C}{\sqrt{L}} \frac{(p - k)^{1-\delta}}{k} \frac{(p - k)^{\delta/2}}{\sqrt{\log(p - k)}},$$

which diverges to infinity, since $k = \mathcal{O}(p^{1-\delta})$.

*Analysis of event* $\mathbb{A}_2$. In order to analyze this event, we first need to condition on the random vector $v := (v^1, \ldots, v^n)$, so as to decouple the random variables $\{D_\ell, \ell \in S\}$. After conditioning on $v$, each variate $nD_\ell, \ell \in S$, is a noncentral $\chi^2_{n,\nu^*}$, with $n$ degrees of freedom and noncentrality parameter $\nu^* = \frac{\beta}{k}\|v\|_2^2$, so that each $D_\ell$ has mean $(\nu^* + n)$.

Since $v$ is a standard Gaussian $n$-vector, we have $\|v\|_2^2 \sim \chi^2_n$. Therefore, if we define the event $\mathbb{B}(v) := \{\frac{\|v\|_2^2}{n} > \frac{3}{2}\}$, the large deviations bound (60a) implies that $\mathbb{P}[\mathbb{B}] \leq \exp(-n/16)$. Therefore, by conditioning on $\mathbb{B}$ and its complement, we obtain

$$\begin{aligned}
(20) \qquad \mathbb{P}[\mathbb{A}_2^c] &\leq \mathbb{P}\left[\min_{\ell \in S} D_\ell > 1 + t \mid \mathbb{B}^c\right] + \mathbb{P}[\mathbb{B}] \\
&\leq (\mathbb{P}[\chi^2_{n,\nu^*} > n(1 + t) \mid \mathbb{B}^c])^k + \exp(-n/16),
\end{aligned}$$

where we have used the conditional independence of $\{D_\ell, \ell \in S\}$. Finally, since $\frac{\|v\|_2^2}{n} \leq \frac{3}{2}$ on the event $\mathbb{B}^c$, we have $\nu^* \leq \frac{3\beta}{2k}n$, and thus

$$\mathbb{P}[\chi^2_{n,\nu^*} > n(1 + t) \mid \mathbb{B}^c] \leq \mathbb{P}\left[\chi^2_{n,\nu^*} > \{n + \nu^*\} + n\left\{t - \frac{3\beta}{2k}\right\} \mid \mathbb{B}^c\right].$$

Since $t = \sqrt{\delta \log(p - k)/n}$ and $n < Lk^2 \log(p - k)$, we have $t \geq \sqrt{\frac{\delta}{L}}\frac{1}{k}$, so that the quantity $\epsilon := \min\{\frac{1}{2}, t - \frac{3\beta}{2k}\}$ is positive for the pre-factor $L > 0$ chosen sufficiently small. Thus, we have

$$\begin{aligned}
\mathbb{P}[\chi^2_{n,\nu^*} > n(1 + t) \mid \mathbb{B}^c] &\leq \mathbb{P}[\chi^2_{n,\nu^*} > \{n + \nu^*\} + n\epsilon] \\
&\leq \exp\left(-\frac{n\epsilon^2}{16(1 + 2(3/2))}\right) = \exp\left(-\frac{n\epsilon^2}{64}\right),
\end{aligned}$$

using the $\chi^2$ tail bound (63). Substituting this upper bound into (20), we obtain

$$\mathbb{P}[\mathbb{A}_2^c] \leq \exp\left(-\frac{kn\epsilon^2}{64}\right) + \exp(-n/16),$$

which certainly vanishes if $\epsilon = \frac{1}{2}$. Otherwise, we have $\epsilon = t - \frac{3\beta}{2k}$ with $t = \sqrt{\frac{\delta \log(p-k)}{n}}$, and we need the quantity

$$\sqrt{kn}\left(t - \frac{3\beta}{2k}\right) = \sqrt{\delta k \log(p - k)} - \frac{3\beta}{2}\sqrt{\frac{n}{k}},$$

to diverge to $+\infty$. This divergence is guaranteed by choosing $n < Lk^2 \log(p - k)$ for $L$ sufficiently small.

**4. Proof of Theorem 2(b).** The proof of our main result is constructive in nature, based on the notion of a *primal–dual certificate*, that is, a primal feasible solution and a dual feasible solution that together satisfy the optimality conditions associated with the SDP (13).

4.1. *High-level proof outline.* We first provide a high-level outline of the main steps in our proof. Under the stated assumptions of Theorem 2, it suffices to construct a rank-one optimal solution $\widehat{Z} = \widehat{z}\,\widehat{z}^T$, constructed from a vector with $\|\widehat{z}\|_2 = 1$, as well as the following properties:

(21a)    Correct sign:    $\operatorname{sign}(\widehat{z}_i) = \operatorname{sign}(z_i^*)$    for all $i \in S$  and

(21b)    Correct exclusion:    $\widehat{z}_j = 0$    for all $j \in S^c$.

Note that our objective function $f(Z) = \operatorname{tr}(\widehat{\Sigma} Z) - \rho_n \sum_{i,j} |Z_{ij}|$ is concave but not differentiable. However, it still possesses a subdifferential (see the books [17, 33] for more details), so that it may be shown that the following conditions are sufficient to verify the optimality of $\widehat{Z} = \widehat{z}\,\widehat{z}^T$.

LEMMA 5. *Suppose that, for each $x \in \mathbb{R}^p$ with $\|x\|_2 = 1$, there exists a sign matrix $\widehat{U} = \widehat{U}(x)$ such that:*

(a) *the matrix $\widehat{U}$ satisfies*

(22)    $$\widehat{U}_{ij} = \begin{cases} \operatorname{sign}(\widehat{z}_i)\operatorname{sign}(\widehat{z}_j), & \text{if } \widehat{z}_i\widehat{z}_j \neq 0, \\ \in [-1, +1], & \text{otherwise;} \end{cases}$$

(b) *the vector $\widehat{z}$ satisfies of $x^T(\widehat{\Sigma} - \rho_n\widehat{U}(x))x \leq \widehat{z}^T(\widehat{\Sigma} - \rho_n\widehat{U}(x))\widehat{z}$.*

*Then $\widehat{Z} = \widehat{z}\widehat{z}^T$ is an optimal rank-one solution.*

PROOF. The subdifferential $\partial f(\widehat{Z})$ of our objective function at $Z = \widehat{Z}$ consists of matrices of the form $\widehat{\Sigma} - \rho_n U$, where $U$ satisfies the condition (22). By the concavity of $f$, for any such $U$ and for all $x \in \mathbb{R}^p$ with $\|x\|_2 = 1$, we have

$$f(xx^T) \leq f(\widehat{Z}) + \operatorname{tr}((\widehat{\Sigma} - \rho_n U)(xx^T - \widehat{Z})).$$

Therefore, it suffices to demonstrate, for each $x \in \mathbb{R}^p$ with $\|x\|_2 = 1$, a valid sign matrix $\widehat{U}(x)$ such that $\operatorname{tr}((\widehat{\Sigma} - \rho_n\widehat{U}(x))(xx^T - \widehat{Z})) \leq 0$. Since we have

$$\operatorname{tr}((\widehat{\Sigma} - \rho_n\widehat{U}(x))xx^T) \leq \operatorname{tr}((\widehat{\Sigma} - \rho_n\widehat{U}(x))\widehat{Z})$$

by assumption (b), the stated conditions are sufficient.  □

REMARKS. Note that if there is a $\widehat{U}$ independent of $x$ such that $\widehat{z}$ satisfies condition (b) of Lemma 5, that is, if $\widehat{z}$ is a maximal eigenvector of $\widehat{\Sigma} - \rho_n \widehat{U}$, then the above argument shows that $\widehat{z} \, \widehat{z}^T$ is in fact "the" optimal solution (i.e., among all matrices in the constraint space, not necessarily rank one).

The condition (22), when combined with the condition (21a), implies that we must have

$$(23) \qquad \widehat{U}_{SS} = \operatorname{sign}(z_S^*) \operatorname{sign}(z_S^*)^T.$$

The remainder of the proof consists in choosing appropriately the remaining dual blocks $\widehat{U}_{SS^c}$ and $\widehat{U}_{S^cS^c}$, and verifying that the primal–dual optimality conditions are satisfied. To describe the remaining steps, it is convenient to define the matrix

$$(24) \qquad \Phi := \widehat{\Sigma} - \rho_n \widehat{U} - \Gamma = \beta z^* z^{*T} - \rho_n \widehat{U} + \Delta,$$

where $\Delta := \widehat{\Sigma} - \Sigma$ is the effective noise in the sample covariance matrix. We divide our proof into three main steps, based on the block structure

$$(25) \quad \Phi = \begin{bmatrix} \Phi_{SS} & \Phi_{SS^c} \\ \Phi_{S^cS} & \Phi_{SS} \end{bmatrix} = \begin{bmatrix} \beta z_S^* z_S^{*T} - \rho_n \widehat{U}_{SS} + \Delta_{SS} & -\rho_n \widehat{U}_{SS^c} + \Delta_{SS^c} \\ -\rho_n \widehat{U}_{S^cS} + \Delta_{S^cS} & -\rho_n \widehat{U}_{S^cS^c} + \Delta_{S^cS^c} \end{bmatrix}.$$

(A) In step A, we analyze the upper-left block $\Phi_{SS}$, using the fixed choice $\widehat{U}_{SS} = \operatorname{sign}(z_S^*) \operatorname{sign}(z_S^*)^T$. We establish conditions on the regularization parameter $\rho_n$ and the noise matrix $\Delta_{SS}$ under which the maximal eigenvector of $\Phi_{SS}$ has the same sign pattern as $z_S^*$. This maximal eigenvector specifies the $k$-dimensional subvector $\widehat{z}_S$ of our optimal primal solution.

(B) In step B, we analyze the off-diagonal block $\Phi_{S^cS}$, in particular establishing conditions on the noise matrix $\Delta_{S^cS}$ under which a valid sign matrix $\widehat{U}_{S^cS}$ can be chosen such that the $p$-vector $\widehat{z} := (\widehat{z}_S, \vec{0}_{S^c})$ is an eigenvector of the full matrix $\Phi$.

(C) In step C, we focus on the lower right block $\Phi_{S^cS^c}$, in particular analyzing conditions on $\Delta_{S^cS^c}$ such that a valid sign matrix $\widehat{U}_{S^cS^c}$ can be chosen such that $\widehat{z}$ defined in step B satisfies condition (b) of Lemma 5.

Our primary interest in this paper is the effective noise matrix $\Delta = \widehat{\Sigma} - \Sigma$ induced by the usual i.i.d. sampling model. However, our results are actually somewhat more general, in that we can provide conditions on arbitrary noise matrices (which need not be of the Wishart type) under which it is possible to construct $(\widehat{z}, \widehat{U})$ as in steps A through C. Accordingly, in order to make the proof as clear as possible, we divide our analysis into two parts: in Section 4.2, we specify sufficient properties on arbitrary noise matrices $\Delta$, and in Section 4.3, we analyze the Wishart ensemble induced by the i.i.d. sampling model and establish sufficient conditions on the sample size $n$. In Section 4.3, we focus exclusively on the special case of the spiked identity covariance, whereas Section 4.4 describes how our results extend to the more general spiked covariance ensembles covered by Theorem 2.

4.2. *Sufficient conditions for general noise matrices.* We now state a series of sufficient conditions, applicable to general noise matrices. So as to clarify the flow of the main proof, we defer the proofs of these technical lemmas to Appendix D.

4.2.1. *Sufficient conditions for step* A. We begin with sufficient condition for the block $(S, S)$. In particular, with the choice (23) of $\widehat{U}_{SS}$ and noting that $\text{sign}(z_S^*) = \sqrt{k} z_S^*$ by assumption, we have

$$\Phi_{SS} = (\beta - \rho_n k) z_S^* z_S^{*T} + \Delta_{SS} := \alpha z_S^* z_S^{*T} + \Delta_{SS},$$

where the quantity $\alpha := \beta - \rho_n k < \beta$ represents a "post-regularization" signal-to-noise ratio. Throughout the remainder of the development, we enforce the constraint

$$(26) \qquad\qquad \rho_n = \frac{\beta}{2k},$$

so that $\alpha = \beta/2$. The following lemma guarantees correct sign recovery [see (21a)], assuming that $\Delta_{SS}$ is "small" in a suitable sense.

LEMMA 6 (Correct sign recovery). *Suppose that the upper-left noise matrix $\Delta_{SS}$ satisfies*

$$(27) \qquad\qquad \|\!|\Delta_{SS}|\!\|_{\infty,\infty} \leq \frac{\alpha}{10} \quad \text{and} \quad \|\!|\Delta_{SS}|\!\|_{2,2} \to 0$$

*with probability 1 as $p \to +\infty$. Then w.a.p. one, the following occurs:*

(a) *The maximal eigenvalue $\gamma_1 := \lambda_{\max}(\Phi_{SS})$ converges to $\alpha$, and its second largest eigenvalue $\gamma_2$ converges to zero.*

(b) *The upper-left block $\Phi_{SS}$ has a unique maximal eigenvector $\widehat{z}_S$ with the correct sign property [i.e., $\text{sign}(\widehat{z}_S) = \text{sign}(z_S^*)$]. More specifically, we have*

$$(28) \qquad\qquad \|\widehat{z}_S - z_S^*\|_\infty \leq \frac{1}{2\sqrt{k}}.$$

4.2.2. *Sufficient conditions for step* B. With the subvector $\widehat{z}_S$ specified, we can now specify the $(p - k) \times k$ submatrix $\widehat{U}_{S^c S}$ so that the vector

$$\widehat{z} := (\widehat{z}_S, \vec{0}_{S^c}) \in \mathbb{R}^p \tag{29}$$

is an eigenvector of the full matrix $\Phi$. In particular, if we define the renormalized quantity $\widetilde{z}_S = \widehat{z}_S / \|\widehat{z}_S\|_1$, and choose

$$\widehat{U}_{S^c S} = \frac{1}{\rho_n} (\Delta_{S^c S} \widetilde{z}_S) \operatorname{sign}(\widehat{z}_S)^T, \tag{30}$$

then some straightforward algebra shows that $(\Delta_{S^c S} - \rho_n \widehat{U}_{S^c S}) \widehat{z}_S = 0$, so that $\widehat{z}$ is an eigenvector of the matrix $\Phi = \beta z^* (z^*)^T - \rho_n \widehat{U} + \Delta$. It remains to verify that the choice (30) is a valid sign matrix (meaning that its entries are bounded in absolute value by one).

LEMMA 7. *Suppose that w.a.p. one, the matrix $\Delta$ satisfies conditions (27), and in addition, for sufficiently small $\delta > 0$, we have*

$$\|\Delta_{S^c S}\|_{\infty, 2} \leq \frac{\delta}{\sqrt{k}}. \tag{31}$$

*Then the specified $\widehat{U}_{S^c S}$ is a valid sign matrix w.a.p. one.*

4.2.3. *Sufficient conditions in step* C. Up to this point, we have established that $\widehat{z} := (\widehat{z}_S, \vec{0}_{S^c})$ is an eigenvector of $\widehat{\Sigma} - \rho_n \widehat{U}$. Thus far, we have specified the sub-blocks $\widehat{U}_{SS}$ and $\widehat{U}_{SS^c}$ of the sign matrix. To complete the proof, it suffices to show that condition (b) in Lemma 5 can be satisfied—namely, that for each $x \in S^{p-1}$, there exists an extension $\widehat{U}_{S^c S^c}(x)$ to our sign matrix such that

$$\widehat{z}^T (\widehat{\Sigma} - \rho_n \widehat{U}(x)) \widehat{z} \geq x^T (\widehat{\Sigma} - \rho_n \widehat{U}(x)) x.$$

Note that it is sufficient to establish the above inequality with $\Phi(x)$ in place of $\widehat{\Sigma} - \rho_n \widehat{U}(x)$.[1] Given any vector $x \in S^{p-1}$, recall the definition (24) of the matrix $\Phi = \Phi(x)$, and observe that $(\widehat{z})^T \Phi(x) \widehat{z} = \gamma_1$ for any choice of $\widehat{U}_{S^c S^c}(x)$. Consider the partition $x = (u, v) \in S^{p-1}$, with $u \in \mathbb{R}^k$ and $v \in \mathbb{R}^m$, where $m = p - k$. We have

$$x^T \Phi x = u^T \Phi_{SS} u + 2 v^T \Phi_{S^c S} u + v^T \Phi_{S^c S^c} v. \tag{32}$$

Let us decompose $u = \mu \widehat{z}_S + \widehat{z}_S^\perp$, where $|\mu| \leq 1$ and $\widehat{z}_S^\perp$ is an element of the orthogonal complement of the span of $\widehat{z}_S$. With this decomposition, we have

$$u^T \Phi_{SS} u = \mu^2 \widehat{z}_S^T \Phi_{SS} \widehat{z}_S + 2 \mu \widehat{z}_S^T \Phi_{SS} \widehat{z}_S^\perp + (\widehat{z}_S^\perp)^T \Phi_{SS} \widehat{z}_S^\perp$$

$$= \mu^2 \gamma_1 + (\widehat{z}_S^\perp)^T \Phi_{SS} \widehat{z}_S^\perp,$$

---

[1] In particular, we have $x^T \Gamma x \leq \|\Gamma\|_{2,2} \|x\|_2^2 = \max\{1, \|\Gamma_{p-k}\|_{2,2}\} \|x\|_2^2 = 1$, while $\widehat{z}^T \Gamma \widehat{z} = \|\widehat{z}_S\|_2^2 = 1$; that is, we have $x^T \Gamma x \leq \widehat{z}^T \Gamma \widehat{z}$.

using the fact that $\widehat{z}_S$ is an eigenvector of $\Phi_{SS}$ with eigenvalue $\gamma_1$ by definition. Note that $\|\widehat{z}_S^\perp\|_2^2 \leq 1 - \mu^2$, so that $(\widehat{z}_S^\perp)^T \Phi_{SS} \widehat{z}_S^\perp$ is bounded by $(1 - \mu^2)\gamma_2$, where $\gamma_2$ is the second largest eigenvalue of $\Phi_{SS}$, which tends to zero according to Lemma 6. We thus conclude that

$$u^T \Phi_{SS} u \leq \mu^2 \gamma_1 + (1 - \mu^2)\gamma_2. \tag{33}$$

The following lemma addresses the remaining two terms in the decomposition (32).

LEMMA 8. *Let* $m = p - k$ *and let* $\mathbb{S} = \{(\eta_i, \ell_i)\}_i$ *be a set of cardinality* $|\mathbb{S}| = \mathcal{O}(m)$. *Suppose that in addition to conditions* (27) *and* (31), *the noise matrix* $\Delta$ *satisfies, w.p. 1,*

$$\max_{\substack{\|v\|_2 \leq \eta, \\ \|v\|_1 \leq \ell}} \sqrt{v^T (\Delta_{S^c S^c} + \Gamma_m) v} \leq \eta + \frac{\delta}{\sqrt{k}}\ell + \varepsilon \qquad \forall (\eta, \ell) \in \mathbb{S}, \tag{34}$$

*for sufficiently small* $\delta, \epsilon > 0$ *as* $m \to +\infty$. *Then w.p. 1, for all* $x \in S^{p-1}$, *there exists a valid sign matrix* $\widehat{U}_{S^c S^c}(x)$ *such that the matrix* $\Phi(x) := \beta z^* z^{*T} - \rho_n \widehat{U}(x) + \Delta$ *satisfies*

$$x^T (\Phi(x)) x \leq \mu^2 \alpha + (1 - \mu^2)\frac{\alpha}{2} \leq \alpha, \tag{35}$$

*where* $|\mu| = |x^T \widehat{z}| \leq 1$.

4.3. *Noise in a sample covariance.* Having established general sufficient conditions on the effective noise matrix, we now turn to the case of i.i.d. samples $x^1, \ldots, x^n$ from the population covariance, and let the effective noise matrix correspond to the difference between the sample and population covariances. Our interest is in providing specific scalings of the triplet $(n, p, k)$ that ensure that the constructions in steps A through C can be carried out. So as to clarify the steps involved, we begin with the proof for the spiked identity ensemble ($\Gamma = I$). In Section 4.4, we provide the extension to nonidentity spiked ensembles.

Recalling our sampling model $x^i = \sqrt{\beta} v^i z^* + g^i$, define the vector $h = \frac{1}{n} \sum_{i=1}^n v^i g^i$. The effective noise matrix $\Delta = \widehat{\Sigma} - \Sigma$ can be decomposed as follows:

$$
\begin{aligned}
\Delta = \beta \underbrace{\left(\frac{1}{n}\sum_{i=1}^n (v^i)^2 - 1\right) z^* z^{*T}}_{P} \\
+ \sqrt{\beta} \underbrace{(z^* h^T + h z^{*T})}_{R} + \underbrace{\left(n^{-1}\sum_{i=1}^n g^i g^{iT} - I_p\right)}_{W}.
\end{aligned}
\tag{36}
$$

We have named each of the three terms that appear in (36), so that we can deal with each one separately in our analysis. The decomposition can be summarized as

$$\Delta = \beta P + \sqrt{\beta} R + W.$$

The last term $W$ is a *centered Wishart random matrix*, whereas the other two are cross terms from the sampling model, involving both random vectors and the unknown eigenvector $z^*$. Defining the standard Gaussian random matrix $G = (g_j^i)_{i,j=1,1}^{n,p} \in \mathbb{R}^{n \times p}$, we can express $W$ concisely as

$$(37) \qquad\qquad W = \frac{1}{n} G^T G - I_p.$$

Our strategy is to examine each of the terms $\beta P$, $\sqrt{\beta} R$ and $W$ separately. For sub-block $\Delta_{SS}$, the corresponding sub-blocks of all the three terms are present, while for sub-block $\Delta_{S^c S}$, only $\sqrt{\beta} R_{S^c S}$ and $W_{S^c S}$ have contributions. Since the conditions to be satisfied by these two sub-blocks are expressed in terms of their (operator) norms, the triangle inequality immediately yields the results for the whole sub-block, once we have established them separately for each of the contributing terms. On the other hand, although the conditions on $\Delta_{S^c S^c}$ (given in Lemma 8) do not have this (sub)additive property, only the Wishart term contributes to this sub-block, and it has a natural decomposition of the form required.

Regarding the Wishart term, the spectral norm ($\|\!|W|\!\|_{2,2}$) of such a random matrix is well characterized [10, 13]; for instance, see claim (38a) in Lemma 10 for one precise statement. The following lemma, concerning the mixed $(\infty, 2)$ norms of submatrices of centered Wishart matrices, is perhaps of independent interest, and plays a key role in our analysis.

LEMMA 9. *Let $W \in \mathbb{R}^{p \times p}$ be a centered Wishart matrix as defined in (37). Let $\mathcal{I}, \mathcal{J} \subset \{1, \ldots, p\}$ be sets of indices, with cardinalities $|\mathcal{I}|, |\mathcal{J}| \to \infty$ as $n, p \to \infty$, and let $W_{\mathcal{I},\mathcal{J}}$ denote the corresponding submatrix. Then, as long as $\max\{|\mathcal{J}|, \log |\mathcal{I}|\}/n = o(1)$, we have*

$$\|\!|W_{\mathcal{I},\mathcal{J}}|\!\|_{\infty,2} = \mathcal{O}\left(\frac{\sqrt{|\mathcal{J}|} + \sqrt{\log |\mathcal{I}|}}{\sqrt{n}}\right)$$

*as $n, p \to +\infty$ with probability 1.*

See Appendix E for the proof of this claim.

4.3.1. *Verifying steps* A *and* B. First, let us look at the Wishart random matrix. The conditions on the upper-left sub-block $W_{SS}$ and lower-left sub-block $W_{S^c S}$ are addressed in the following lemma.

Lemma 10. *As $(n, p, k) \to +\infty$, we have w.a.p. one*

$$\|\|W_{SS}\|\|_{2,2} = \mathcal{O}\left(\sqrt{\frac{k}{n}}\right), \tag{38a}$$

$$\|\|W_{SS}\|\|_{\infty,\infty} = \mathcal{O}\left(\sqrt{\frac{k^2}{n}}\right), \tag{38b}$$

$$\|\|W_{S^cS}\|\|_{\infty,2} = \mathcal{O}\left(\frac{\sqrt{k} + \sqrt{\log(p-k)}}{\sqrt{n}}\right). \tag{38c}$$

*In particular, under the scaling $n > Lk \log(p-k)$ and $k = \mathcal{O}(\log p)$, the conditions of Lemmas 6 and 7 are satisfied for $W_{SS}$ and $W_{S^cS}$ for sufficiently large $L$.*

Proof. Assertion (38a) about the spectral norm of $W_{SS}$ follows directly from known results on singular values of Gaussian random matrices (e.g., see [10, 13]). To bound the mixed norm $\|\|W_{S^cS}\|\|_{\infty,2}$, we apply Lemma 9 with the choices $\mathcal{I} = S^c$ and $\mathcal{J} = S$, noting that $|\mathcal{I}| = p - k$ and $|\mathcal{J}| = k$. Finally, to obtain a bound on $\|\|W_{SS}\|\|_{\infty,\infty}$, we first bound $\|\|W_{SS}\|\|_{\infty,2}$. Again using Lemma 9, this time with the choices $\mathcal{I} = \mathcal{J} = S$, we obtain

$$\|\|W_{SS}\|\|_{\infty,2} = \mathcal{O}\left(\frac{\sqrt{k} + \sqrt{\log k}}{\sqrt{n}}\right) = \mathcal{O}\left(\sqrt{\frac{k}{n}}\right) \tag{39}$$

as $n, k \to \infty$. Now, using the fact that for any $x \in \mathbb{R}^k$, $\|x\|_2 \le \sqrt{k}\|x\|_\infty$, we obtain

$$\|\|W_{SS}\|\|_{\infty,\infty} = \max_{\|x\|_\infty \le 1} \|W_{SS}x\|_\infty \le \max_{\|x\|_2 \le \sqrt{k}} \|W_{SS}x\|_\infty = \sqrt{k}\|\|W_{SS}\|\|_{\infty,2}.$$

Combined with the inequality (39), we obtain the stated claim (38b). □

We now turn to the cross-term $R$, and establish the following result.

Lemma 11. *The matrix $R = z^*h^T + hz^{*T}$, as defined in (36), satisfies the conditions of Lemmas 6 and 7.*

Proof. First observe that $h$ may be viewed as a vector consisting of the off-diagonal elements of the first column of a $(p+1) \times (p+1)$ Wishart matrix, say $W'$. This representation follows since $h_j = \frac{1}{n}\sum_{i=1}^n v^i g_j^i$, where the Gaussian variable $v^i$ is independent of $g_j^i$ for all $1 \le j \le p$. For ease of reference, let us index rows and columns of $W'$ by $1', 1, \ldots, p$, let $S' = \{1'\} \cup S$, and let $h = W'_{1', S \cup S^c}$. (Recall that $S \cup S^c$ is simply $\{1, \ldots, p\}$.)

Since the spectral norm of a matrix is an upper bound on the $\ell_2$-norm of any column, we have

$$\|h_S\|_2 \leq \|W'_{S'S'}\|_{2,2} = \mathcal{O}\left(\sqrt{\frac{k+1}{p}}\right),\tag{40}$$

where we used known bounds [10] on singular values of Gaussian random matrices. Under the scaling $n > Lk\log(p-k)$, we thus have $\|h_S\|_2 \xrightarrow{P} 0$. By Lemma 15, we have $\mathbb{P}[|W'_{ij}| > t] \leq C\exp(-cnt^2)$ for $t > 0$ sufficiently small, which implies (via union bound) that

$$\|h\|_\infty = \mathcal{O}\left(\sqrt{\frac{\log(p)}{n}}\right) = \mathcal{O}\left(\frac{1}{\sqrt{k}}\right),\tag{41}$$

under our assumed scaling. Note also that $\|h\|_\infty = \max\{\|h_S\|_\infty, \|h_{S^c}\|_\infty\}$, that is, the $\infty$-norm of each of these subvectors are also $\mathcal{O}(k^{-1/2})$. Assume for the following that $L$ is chosen large enough so that $\|h\|_\infty \leq \delta/\sqrt{k}$.

Now, to complete the proof, let us first examine the spectral norm of $R_{SS} = z_S^* h_S^T + h_S z_S^{*T}$. The two (possibly) nonzero eigenvalues of this matrix are $z_S^{*T} h_S \pm \|z_S^*\|_2\|h_S\|_2$, whence we have

$$\|R_{SS}\|_{2,2} \leq |z_S^{*T} h_S| + \|z_S^*\|_2\|h_S\|_2 \leq 2\|h_S\|_2 \xrightarrow{P} 0.$$

As for the (matrix) $\infty$-norm of $R_{SS}$, let us exploit the "maximum row sum" interpretation, that is, $\|R_{SS}\|_{\infty,\infty} = \max_{i\in S}\sum_{j\in S}|R_{ij}|$ (cf. Appendix A) to deduce

$$\begin{aligned}
\|R_{SS}\|_{\infty,\infty} &\leq \|z_S^* h_S^T\|_{\infty,\infty} + \|h_S z_S^{*T}\|_{\infty,\infty}\\
&\leq \left(\max_{i\in S}|z_i^*|\right)\|h_S^T\|_1 + \left(\max_{i\in S}|h_i|\right)\|z_S^{*T}\|_1\\
&\leq \frac{1}{\sqrt{k}}\|W'_{S'S'}\|_{\infty,\infty} + \|h_S\|_\infty\sqrt{k}.
\end{aligned}$$

From the argument of Lemma 10, we have $\|W'_{S'S'}\|_{\infty,\infty} = \mathcal{O}(\sqrt{\frac{k^2}{n}})$, so that

$$\frac{1}{\sqrt{k}}\|W'_{S'S'}\|_{\infty,\infty} = \mathcal{O}\left(\sqrt{\frac{k}{n}}\right) \xrightarrow{P} 0$$

and moreover, the norm $\|R_{SS}\|_{\infty,\infty}$ can be made smaller than $2\delta$, by choosing $L$ sufficiently large in the relation $n > Lk\log(p-k)$.

Finally, to establish the additional condition required by Lemma 7—namely (31)—notice that

$$\|R_{S^cS}\|_{\infty,2} = \max_{\|y\|_2=1}\|R_{S^cS}y\|_\infty$$

$$= \max_{\|y\|_2=1} \|h_{S^c} z_S^{*T} y\|_\infty$$

$$= \Big( \max_{\|y\|_2=1} |z_S^{*T} y| \Big) \|h_{S^c}\|_\infty \le \frac{\delta}{\sqrt{k}},$$

where the last line uses $\max_{\|y\|_2=1} |z_S^{*T} y| = \|z_S^*\|_2 = 1$, thereby completing the proof. $\square$

Finally, we examine the first term in (36), that is, $P$. As this term only contributes to the upper-left block, we only need to establish that it satisfies Lemma 6.

LEMMA 12. *The matrix $P_{SS}$ satisfies condition (27) of Lemma 6.*

PROOF. Note that for any matrix norm, we have $\|\|P_{SS}\|\| = |n^{-1} \sum_{i=1}^n (v^i)^2 - 1| \|\|z_S^* z_S^{*T}\|\|$. Now, notice that $\|\|z_S^* z_S^{*T}\|\|_{2,2} = |z_S^{*T} z_S^*| = 1$. Also, using the "maximum row sum" characterization of matrix $\infty$-norm, we have $\|\|z_S^* z_S^{*T}\|\|_{\infty,\infty} = \sum_{j=1}^k |(\pm\frac{1}{\sqrt{k}})(\pm\frac{1}{\sqrt{k}})| = 1$. Now by the strong law of large numbers, $|n^{-1} \times \sum_{i=1}^n (v^i)^2 - 1| \overset{\text{a.s.}}{\to} 0$ as $n \to \infty$. It follows that with probability 1

$$\|\|P_{SS}\|\|_{2,2} = \|\|P_{SS}\|\|_{\infty,\infty} \to 0,$$

which clearly implies condition (27). $\square$

4.3.2. *Verifying step* C. For this step, we only need to consider the lower-right block of $W$; that is, we only need to verify condition (34) of Lemma 8 for $\Delta_{S^c S^c} = W_{S^c S^c}$. Recall that $W = n^{-1} G^T G - I_p$ where $G$ is a $n \times p$ (canonical) Gaussian matrix [see (37)]. With a slight abuse of notation, let $G_{S^c} = (G_{ij})$ for $1 \le i \le n$ and $j \in S^c$. Note that $G_{S^c} \in \mathbb{R}^{n \times m}$ where $m = p - k$ and

$$\Delta_{S^c S^c} + I_m = W_{S^c S^c} + I_m = n^{-1} G_{S^c}^T G_{S^c}.$$

Now, we can simplify the quadratic form in (34) as

$$\sqrt{v^T(\Delta_{S^c S^c} + I_m)v} = \sqrt{\|n^{-1/2} G_{S^c} v\|_2^2} = \|n^{-1/2} G_{S^c} v\|_2$$

for which we have the following lemma.

LEMMA 13. *For any $M > 0$ and $\varepsilon > 0$, there exists a constant $B > 0$ such that for any set $\mathbb{S} = \{(\eta_i, \ell_i)\}_i$ with elements in $(0, M) \times \mathbb{R}^+$ and cardinality $|\mathbb{S}| = \mathcal{O}(m)$, we have*

$$(42) \qquad \max_{\substack{\|v\|_2 \le \eta, \\ \|v\|_1 \le \ell}} \|n^{-1/2} G_{S^c} v\|_2 \le \eta + B\sqrt{\frac{\log m}{n}} \ell + \varepsilon \qquad \forall (\eta, \ell) \in \mathbb{S},$$

*as $p \to \infty$, with probability 1. In particular, under the scaling $n > Lk \log m$, condition (34) of Lemma 8 is satisfied for $L$ large enough.*

PROOF. Without loss of generality, assume $M = 1$. We begin by controlling the expectation of the left-hand side, using an argument based on the Gordon–Slepian theorem [26], similar to that used for establishing bounds on spectral norms of random Gaussian matrices (e.g., [10]). First, we require some notation: for a zero-mean random variable $Z$, define its standard deviation $\sigma(Z) = (\mathbb{E}|Z|^2)^{1/2}$. For vectors $x, y$ of the same dimension, define the Euclidean inner product $\langle x, y \rangle = x^T y$. For matrices $X, Y$ of the same dimension (although not necessarily symmetric), recall the Hilbert–Schmidt norm

$$\|X\|_{\mathrm{HS}} := \langle\!\langle X, X \rangle\!\rangle^{1/2} = \left( \sum_{i,j} X_{ij}^2 \right)^{1/2}.$$

Given some (possibly uncountable) index set $\{t \in T\}$, let $(X_t)_{t \in T}$ and $(Y_t)_{t \in T}$ be a pair of centered Gaussian processes. One version of the Gordon–Slepian theorem (see [26]) asserts that if $\sigma(X_s - X_t) \le \sigma(Y_s - Y_t)$ for all $s, t \in T$, then we have

$$(43) \qquad \mathbb{E}\Big[\sup_{t \in T} X_t\Big] \le \mathbb{E}\Big[\sup_{t \in T} Y_t\Big].$$

For simplicity in notation, define $\widetilde{H} := G_{S^c} \in \mathbb{R}^{n \times m}$, $H := n^{-1/2} G_S^c$, and fix some $\eta, \ell > 0$. We wish to bound

$$f(\widetilde{H}; \eta, \ell) := \max_{\substack{\|v\|_2 \le \eta, \\ \|v\|_1 \le \ell}} \|\widetilde{H}v\|_2 = \max_{\substack{\|v\|_2 \le \eta, \\ \|v\|_1 \le \ell, \\ \|u\|_2 = 1}} \langle \widetilde{H}v, u \rangle,$$

where $v \in \mathbb{R}^m$, $u \in \mathbb{R}^n$. Note that $\langle \widetilde{H}v, u \rangle = u^T \widetilde{H}v = \mathrm{tr}(\widetilde{H}vu^T) = \langle\!\langle \widetilde{H}, uv^T \rangle\!\rangle$. Consider $\widetilde{H}$ to be a (canonical) Gaussian vector in $\mathbb{R}^{mn}$, take

$$(44) \qquad T := \{t = (u, v) \in \mathbb{R}^n \times \mathbb{R}^m \mid \|v\|_2 \le \eta, \|v\|_1 \le \ell, \|u\|_2 = 1\}$$

and define $X_t = \langle\!\langle \widetilde{H}, uv^T \rangle\!\rangle$ for $t \in T$. Observe that $(X_t)_{t \in T}$ is a (centered) canonical Gaussian process generated by $\widetilde{H}$, and $f(\widetilde{H}; \eta, \ell) = \max_{t \in T} X_t$. We compare this to the maximum of another Gaussian process $(Y_t)_{t \in T}$, defined as $Y_t = \langle (g, h), (u, v) \rangle$ where $g \in \mathbb{R}^n$ and $h \in \mathbb{R}^m$ are Gaussian vectors with $\mathbb{E}[gg^T] = \eta^2 I_n$ and $\mathbb{E}[hh^T] = I_m$. Note that, for example,

$$\sigma(\langle g, u \rangle) = (\mathbb{E}\langle g, u \rangle^2)^{1/2} = (u^T \mathbb{E}[gg^T]u)^{1/2} = \eta \|u\|_2,$$

in which the left-hand size is the norm of a process $(\langle g, u \rangle)_u$ expressed in terms of the norm of a vector (i.e., its index).

Let $t = (u, v) \in T$ and $t' = (u', v') \in T$. Assume, without loss of generality, that $\|v'\|_2 \le \|v\|_2$. Then, we have

$$
\begin{aligned}
\sigma^2(X_t - X_{t'}) &= \|uv^T - u'v'^T\|_{\mathrm{HS}}^2 \\
&= \|uv^T - u'v^T + u'v^T - u'v'^T\|_{\mathrm{HS}}^2 \\
&= \|v\|_2^2\|u - u'\|_2^2 + \|u'\|_2^2\|v - v'\|_2^2 \\
&\quad + 2(u^Tu' - \|u'\|_2^2)(\|v\|_2^2 - v^Tv') \\
&\le \eta^2\|u - u'\|_2^2 + \|v - v'\|_2^2 = \sigma^2(Y_t - Y_{t'}),
\end{aligned}
$$

where we have used Cauchy–Schwarz inequality to deduce $|u^Tu'| \le 1 = \|u'\|_2^2$ and $|v^Tv'| \le \|v\|_2\|v'\|_2 \le \|v\|_2^2$. Thus, the Gordon–Slepian lemma is applicable, and we obtain

$$
\begin{aligned}
\mathbb{E}f(\widetilde{H}; \eta, \ell) &\le \mathbb{E}\max_{t \in T} Y_t \\
&= \mathbb{E}\max_{\|u\|_2 = 1}\langle g, u\rangle + \mathbb{E}\max_{\substack{\|v\|_2 \le \eta, \\ \|v\|_1 \le \ell}}\langle h, v\rangle \\
&\le \mathbb{E}\|g\|_2 + (\mathbb{E}\|h\|_\infty)\ell \\
&< \sqrt{n}\eta + (\sqrt{3\log m})\ell,
\end{aligned}
$$

where we have used $(\mathbb{E}\|g\|_2)^2 < \mathbb{E}(\|g\|_2^2) = \mathbb{E}\operatorname{tr}(gg^T) = \operatorname{tr}\mathbb{E}(gg^T) = n\eta^2$; the bound used for $\mathbb{E}\|h\|_\infty$ follows from standard Gaussian tail bounds [26]. Noting that $H = n^{-1/2}\widetilde{H}$, we obtain $\mathbb{E}f(H; \eta, \ell) \le \eta + \sqrt{\frac{3\log m}{n}}\ell$.

The final step is to argue that $f(H; \eta, \ell)$ is sufficiently close to its mean. For this, we will use concentration of Gaussian measure [25, 26] for Lipschitz functions in $\mathbb{R}^{mn}$. To see that $A \to f(A; \eta, \ell)$ is in fact 1-Lipschitz, note that it satisfies the triangle inequality and it is bounded above by the spectral norm. Thus,

$$
|f(\widetilde{H}; \eta, \ell) - f(\widetilde{F}; \eta, \ell)| \le f(\widetilde{H} - \widetilde{F}; \eta, \ell) \le \|\widetilde{H} - \widetilde{F}\|_{2,2} \le \|\widetilde{H} - \widetilde{F}\|_{\mathrm{HS}},
$$

where we have used the assumption $\eta \le 1$. Noting that $H = n^{-1/2}\widetilde{H}$ and $f(H; \eta, \ell) = n^{-1/2}f(\widetilde{H}; \eta, \ell)$, Gaussian concentration of measure for 1-Lipschitz functions [25] implies that

$$
\mathbb{P}[f(H; \eta, \ell) - \mathbb{E}[f(H; \eta, \ell)] > t] \le \exp(-nt^2/2).
$$

Finally, we use union bound to establish the result uniformly over $\mathbb{S}$. By assumption, there exists some $K > 0$ such that $|\mathbb{S}| \le Km$. Thus,

$$
\mathbb{P}\Big[\max_{(\eta, \ell) \in \mathbb{S}}(f(H; \eta, \ell) - (\eta + \sqrt{(3\log m)/n}\cdot\ell)) > t\Big] \le K\exp(-nt^2/2 + \log m).
$$

Now, fix some $\varepsilon > 0$, take $t = \sqrt{\frac{6 \log m}{n}}$ and apply the Borell–Cantelli lemma to conclude that

$$\max_{(\eta, \ell) \in \mathbb{S}} \left[ f(H; \eta, \ell) - \left( \eta + \sqrt{\frac{3 \log m}{n}} \cdot \ell \right) \right] \leq \sqrt{\frac{6 \log m}{n}} \leq \varepsilon,$$

eventually (w.p. 1). $\square$

### 4.4. *Nonidentity noise covariance.*

In this section, we specify how the proof is extended to (population) covariance matrices having a more general base covariance term $\Gamma_{p-k}$ in (5). Let $\Gamma_{p-k}^{1/2}$ denote the (symmetric) square root of $\Gamma_{p-k}$. We can write samples from this model as

$$(45) \qquad \tilde{x}^i = \sqrt{\beta} v^i z^* + \tilde{g}^i, \qquad i = 1, \ldots, n,$$

where

$$(46) \qquad \tilde{g}^i = \begin{pmatrix} g_S^i \\ \Gamma_{p-k}^{1/2} g_{S^c}^i \end{pmatrix}$$

with $g^i \sim N(0, I_p)$ and $v^i \sim N(0, 1)$ standard independent Gaussian random variables.

Denoting the resulting sample covariance as $\widehat{\Sigma}$, we can obtain an expression for the noise matrix $\Delta = \widehat{\Sigma} - \Sigma$. The result will be similar to expansion (36) with $h$ and $W$ appropriately modified; more specifically, we have

$$(47) \qquad \tilde{h}_S = h_S, \qquad \tilde{h}_{S^c} = \Gamma_{p-k}^{1/2} h_{S^c},$$

$$(48) \qquad \widetilde{W}_{SS} = W_{SS}, \qquad \widetilde{W}_{S^c S} = \Gamma_{p-k}^{1/2} W_{S^c S}, \qquad \widetilde{W}_{S^c S^c} = \Gamma_{p-k}^{1/2} W_{S^c S^c} \Gamma_{p-k}^{1/2}.$$

Note that the $P$-term is unaffected.

Re-examining the proof presented for the case $\Gamma_{p-k} = I_{p-k}$, we can identify conditions imposed on $h$ and $W$ to guarantee optimality. By imposing sufficient constraints on $\Gamma_{p-k}$, we can make $\tilde{h}$ and $\widetilde{W}$ satisfy the same conditions. The rest of the proof will then be exactly the same as the case $\Gamma_{p-k} = I_{p-k}$. As before, we proceed by verifying steps A through C in sequence.

#### 4.4.1. *Verifying steps* A *and* B.

Examining the proof of Lemma 11, we observe that we need bounds on $\|\tilde{h}_S\|_2, \|\tilde{h}_S\|_1$ and $\|\tilde{h}\|_\infty = \max\{\|\tilde{h}_S\|_\infty, \|\tilde{h}_{S^c}\|_\infty\}$. Since $\tilde{h}_S = h_S$, we should only be concerned with $\|\tilde{h}_{S^c}\|_\infty$, for which we simply have

$$\|\tilde{h}_{S^c}\|_\infty \leq \|\Gamma_{p-k}^{1/2}\|_{\infty,\infty} \|h_{S^c}\|_\infty.$$

Thus, assumption (6a)—that is, $\|\Gamma^{1/2}\|_{\infty,\infty} = \mathcal{O}(1)$—guarantees that Lemma 11 also holds for (nonidentity) $\Gamma$.

Similarly, for Lemma 10 to hold, we need to investigate $\||\widetilde{W}_{S^c S}\||_{\infty,2}$, since this is the only norm (among those considered in the lemma) affected by a nonidentity $\Gamma$. Using sub-multiplicative property of operator norms [see relation (58) in Appendix A], we have

$$\||\widetilde{W}_{S^c S}\||_{\infty,2} \leq \||\Gamma_{p-k}^{1/2}\||_{\infty,\infty} \||W_{S^c S}\||_{\infty,2},$$

so that the same boundedness assumption (6a) is sufficient.

4.4.2. *Verifying step* C. For the lower-right block $\widetilde{W}_{S^c S^c}$, we first have to verify Lemma 13. We also need to examine the proof of Lemma 8 where the result of Lemma 13—namely relation (42)—was used. Let $\widetilde{G} = (\tilde{g}_j^i)_{i,j=1,1}^{n,p}$ and let $\widetilde{G}_{S^c} = (\widetilde{G}_{ij})$ for $1 \leq i \leq n$ and $j \in S^c$. Note that $\widetilde{G}_{S^c}^T \in \mathbb{R}^{(p-k)\times n}$ and we have

$$\widetilde{G}_{S^c}^T = (\tilde{g}_{S^c}^1, \ldots, \tilde{g}_{S^c}^n) = \Gamma_{p-k}^{1/2}(g_{S^c}^1, \ldots, g_{S^c}^n) = \Gamma_{p-k}^{1/2} G_{S^c}^T.$$

Using this notation, we can write $\widetilde{W}_{S^c S^c} = n^{-1}\widetilde{G}_{S^c} - \Gamma_{p-k} = \Gamma_{p-k}^{1/2}(n^{-1}G_{S^c}^T G_{S^c} - I_{p-k})\Gamma_{p-k}^{1/2}$, consistent with (48).

Now to establish a version of (42), we have to consider the maximum of

$$\|n^{-1/2}\widetilde{G}_{S^c} v\|_2 = \|n^{-1/2}G_{S^c}\Gamma_{p-k}^{1/2} v\|_2$$

over the set where $\|v\|_2 \leq \eta$ and $\|v\|_1 \leq \ell$. Let $\tilde{v} = \Gamma_{p-k}^{1/2} v$ and note that for any consistent pair of vector–matrix norms we have $\|\tilde{v}\| \leq \||\Gamma_{p-k}^{1/2}\|| \|v\|$. Thus, for example, $\|v\|_2 \leq \eta$ implies $\|\tilde{v}\|_2 \leq \||\Gamma_{p-k}^{1/2}\||_{2,2}\eta$, and similarly for the $\ell_1$-norm. Now, if we assume that Lemma 13 holds for $G_{S^c}$, we obtain, for all $(\eta, \ell) \in \mathbb{S}$, the inequality

$$
\begin{aligned}
\max_{\substack{\|v\|_2 \leq \eta, \\ \|v\|_1 \leq \ell}} \|n^{-1/2}\widetilde{G}_{S^c} v\|_2 &\leq \max_{\substack{\|\tilde{v}\|_2 \leq \||\Gamma_{p-k}^{1/2}\||_{2,2}\eta, \\ \|\tilde{v}\|_1 \leq \||\Gamma_{p-k}^{1/2}\||_{1,1}\ell}} \|n^{-1/2}G_{S^c}\tilde{v}\|_2 \\
&\leq \||\Gamma_{p-k}^{1/2}\||_{2,2}\eta + B\||\Gamma_{p-k}^{1/2}\||_{1,1}\sqrt{\frac{\log m}{n}}\ell + \varepsilon.
\end{aligned}
\tag{49}
$$

Thus, one observes that the boundedness condition (6a) guarantees that

$$\||\Gamma_{p-k}^{1/2}\||_{1,1} = \||\Gamma_{p-k}^{1/2}\||_{\infty,\infty} \leq A_1,$$

thereby taking care of the second term in (49). More specifically, the constant $A_1$ is simply absorbed into some $B' = BA_1$. In addition, we also require a bound on $\||\Gamma_{p-k}^{1/2}\||_{2,2}$, which follows from our assumption $\||\Gamma_{p-k}\||_{2,2} \leq 1$.

However, the fact that the factor multiplying $\eta$ in (49) is no longer unity has to be addressed more carefully.

Recall that inequality (42) was used in the proof of Lemma 8 to establish a bound on

$$v^{*T}\Delta_{S^cS^c}v^* = v^{*T}W_{S^cS^c}v^* = v^{*T}(H^TH - I_{p-k})v^* = \|Hv^*\|_2^2 - \|v^*\|_2^2,$$

where $H = n^{-1/2}G_{S^c}$. The bound obtained on this term is given by (76). We focus on the core idea, omitting some technical details such as the discretization argument.[2] Replacing $W_{S^cS^c}$ with $\widetilde{W}_{S^cS^c}$, we need to establish a similar bound on

$$v^{*T}\widetilde{W}_{S^cS^c}v^* = v^{*T}(n^{-1}\widetilde{G}_{S^c}^T\widetilde{G}_{S^c} - \Gamma_{p-k})v^* = \|n^{-1/2}\widetilde{G}_{S^c}v^*\|_2^2 - \|\Gamma_{p-k}^{1/2}v^*\|_2^2.$$

Note that $\|v^*\|_2 \leq \|\Gamma_{p-k}^{-1/2}\|_{2,2}\|\Gamma_{p-k}^{1/2}v^*\|_2$ or, equivalently, $\|\Gamma_{p-k}^{-1/2}\|_{2,2}^{-1}\|v^*\|_2 \leq \|\Gamma_{p-k}^{1/2}v^*\|_2$. Thus, using (49), one obtains

$$\|n^{-1/2}\widetilde{G}_{S^c}v^*\|_2^2 - \|\Gamma_{p-k}^{1/2}v^*\|_2^2 \leq (\|\Gamma_{p-k}^{1/2}\|_{2,2}^2 - \|\Gamma_{p-k}^{-1/2}\|_{2,2}^{-2})\|v^*\|_2^2$$
$$+ \text{(terms of lower order in } \|v^*\|_2).$$

Note that unlike the case $\Gamma_{p-k} = I_{p-k}$, the term quadratic in $\|v^*\|_2$ does not vanish in general. Thus, we have to assume that its coefficient is eventually small compared to $\beta$. More specifically, we assume

$$(50) \qquad \|\Gamma_{p-k}^{1/2}\|_{2,2}^2 - \|\Gamma_{p-k}^{-1/2}\|_{2,2}^{-2} \leq \frac{\alpha}{4} = \frac{\beta}{8}, \qquad \text{eventually.}$$

The boundedness assumptions on $\|\Gamma_{p-k}^{1/2}\|_{1,1}$ and $\|\Gamma_{p-k}^{1/2}\|_{2,2}$ now allows for the rest of the terms to be made less than $\alpha/4$, using arguments similar to the proof of Lemma 8, so that the overall objective is less than $\alpha/2$, eventually. This concludes the proof.

Noting that $\|\Gamma_{p-k}^{1/2}\|_{2,2}^2 = \lambda_{\max}(\Gamma_{p-k})$ and $\|\Gamma_{p-k}^{-1/2}\|_{2,2}^{-2} = \lambda_{\min}(\Gamma_{p-k})$, we can summarize the conditions sufficient for Lemma 8 to extend to general covariance structure as follows:

$$(51a) \qquad\qquad\qquad \|\Gamma_{p-k}^{1/2}\|_{1,1} = \|\Gamma_{p-k}^{1/2}\|_{\infty,\infty} = \mathcal{O}(1);$$

$$(51b) \qquad\qquad\qquad \lambda_{\max}(\Gamma_{p-k}) \leq 1;$$

$$(51c) \qquad\qquad \lambda_{\max}(\Gamma_{p-k}) - \lambda_{\min}(\Gamma_{p-k}) \leq \frac{\beta}{8}$$

as stated previously.

---

[2]In particular, we will assume that $v^*$ saturates (49), so that $\|v^*\|_2 = \eta$. For a more careful argument see the proof of Lemma 8.

**5. Proof of Theorem 3.** Our proof is based on the standard approach of applying Fano's inequality (e.g., [7, 16, 37, 38]). Let $\mathbb{S}$ denote the collection of all possible support sets, that is, the collection of $k$-subsets of $\{1, \ldots, p\}$ with cardinality $|\mathbb{S}| = \binom{p}{k}$; we view $S$ as a random variable distributed uniformly over $\mathbb{S}$. Let $\mathbb{P}_S$ denote the distribution of a sample $X \sim N(0, \Sigma_p(S))$ from a spiked covariance model, conditioned on the maximal eigenvector having support set $S$, and let $X^n = (x^1, \ldots, x^n)$ be a set of $n$ i.i.d. samples. In information-theoretic terms, we view any method of support recovery as a decoder that operates on the data $X^n$ and outputs an estimate of the support $\widehat{S} = \phi(X^n)$—in short, a (possibly random) map $\phi : (\mathbb{R}^p)^n \to \mathbb{S}$. Using the 0–1 loss to compare an estimate $\widehat{S}$ and the true support set $S$, the associated risk is simply the probability of error $\mathbb{P}[\text{error}] = \sum_{S \in \mathbb{S}} \frac{1}{\binom{p}{k}} \mathbb{P}_S[\widehat{S} \neq S]$. Due to symmetry of the ensemble, in fact we have $\mathbb{P}[\text{error}] = \mathbb{P}_S[\widehat{S} \neq S]$, where $S$ is some fixed but arbitrary support set, a property that we refer to as *risk flatness*.

In order to generate suitably tight lower bounds, we restrict attention to the following *sub-collection* $\widetilde{\mathbb{S}}$ of support sets:

$$\widetilde{\mathbb{S}} := \{S \in \mathbb{S} \mid \{1, \ldots, k-1\} \subset S\},$$

consisting of those $k$-element subsets that contain $\{1, \ldots, k-1\}$ and one element from $\{k, \ldots, p\}$. By risk flatness, the probability of error with $S$ chosen uniformly at random from the original ensemble $\mathbb{S}$ is the same as the probability of error with $S$ chosen uniformly from $\widetilde{\mathbb{S}}$. Letting $U$ denote a subset chosen uniformly at random from $\widetilde{\mathbb{S}}$, using Fano's inequality, we have the lower bound

$$\mathbb{P}[\text{error}] \geq 1 - \frac{I(U; X^n) + \log 2}{\log |\widetilde{\mathbb{S}}|},$$

where $I(U; X^n)$ is the mutual information between the data $X^n$ and the randomly chosen support set $U$, and $|\widetilde{\mathbb{S}}| = p - k + 1$ is the cardinality of $\widetilde{\mathbb{S}}$.

It remains to obtain an upper bound on $I(U; X^n) = H(X^n) - H(X^n|U)$. By chain rule for entropy, we have $H(X^n) \leq nH(x)$. Next, using the maximum entropy property of the Gaussian distribution [7], we have

$$(52) \qquad H(X^n) \leq nH(x) \leq n \left\{ \frac{p}{2}[1 + \log(2\pi)] + \frac{1}{2} \log \det \mathbb{E}[xx^T] \right\},$$

where $\mathbb{E}[xx^T]$ is the covariance matrix of $x$. On the other hand, given $U = \overline{U}$, the vector $X^n$ is a collection of $n$ Gaussian $p$-vectors with covariance matrix $\Sigma_p(\overline{U})$. The determinant of this matrix is $1 + \beta$, independent of $\overline{U}$, so that we have

$$(53) \qquad H(X^n|U) = \frac{np}{2}[1 + \log(2\pi)] + \frac{n}{2} \log(1 + \beta).$$

Combining (52) and (53), we obtain

$$(54) \qquad I(U; X^n) \leq \frac{n}{2}\{\log \det \mathbb{E}[xx^T] - \log(1 + \beta)\}.$$

The following lemma, proved in Appendix F, specifies the form of the log determinant of the covariance matrix $\Sigma_M := \mathbb{E}[xx^T]$.

LEMMA 14. *The log determinant has the exact expression*

$$
\begin{aligned}
(55) \qquad \log \det \Sigma_M = {} & \log(1 + \beta) + \log\left(1 - \frac{\beta}{1 + \beta}\frac{p - k}{k(p - k + 1)}\right) \\
& + (p - k)\log\left(1 + \frac{\beta}{k(p - k + 1)}\right).
\end{aligned}
$$

Substituting (55) into (54) and using the inequality $\log(1 + \alpha) \leq \alpha$, we obtain

$$
\begin{aligned}
I(U; & X^n) \\
& \leq \frac{n}{2}\left\{\log\left(1 - \frac{\beta}{1 + \beta}\frac{p - k}{k(p - k + 1)}\right) + (p - k)\log\left(1 + \frac{\beta}{k(p - k + 1)}\right)\right\} \\
& \leq \frac{n}{2}\left\{-\frac{\beta}{1 + \beta}\frac{p - k}{k(p - k + 1)} + \frac{\beta(p - k)}{k(p - k + 1)}\right\} \\
& = \frac{n}{2}\left\{\frac{\beta^2}{1 + \beta}\frac{p - k}{k(p - k + 1)}\right\} \\
& \leq \frac{\beta^2}{2(1 + \beta)}\frac{n}{k}.
\end{aligned}
$$

From the Fano bound (52), the error probability is greater than $\frac{1}{2}$ if $\frac{\beta^2}{1 + \beta}\frac{n}{k} < \log(p - k) < \log |\widetilde{\mathbb{S}}|$, which completes the proof.

**6. Discussion.** In this paper, we studied the problem of recovering the support of a sparse eigenvector in a spiked covariance model. Our analysis allowed for high-dimensional scaling, where the problem size $p$ and sparsity index $k$ increase as functions of the sample size $n$. We analyzed two computationally tractable methods for sparse eigenvector recovery—diagonal thresholding and a semidefinite programming (SDP) relaxation [9]—and provided precise conditions on the scaling of the triplet $(n, p, k)$ under which they succeed (or fail) in correctly recovering the support. The probability of success using diagonal thresholding undergoes a phase transition in terms of the rescaled sample size $\theta_{\mathrm{dia}}(n, p, k) = n/(k^2 \log(p - k))$, whereas the more complex SDP relaxation, when it has a rank-one solution, succeeds once the rescaled sample size $\theta_{\mathrm{sdp}}(n, p, k) = n/(k \log(p - k))$ is sufficiently large. Thus, the SDP relaxation has greater statistical efficiency, by a factor of $k$ relative to the simple diagonal thresholding method, but also a substantially larger computational complexity. Finally, using information-theoretic methods, we showed that no method, regardless of its computational complexity, can recover the support set with vanishing error probability if $\theta_{\mathrm{sdp}}(n, p, k)$ is smaller than a critical constant. Our results thus provide some insight into the trade-offs between statistical and computational efficiency in high-dimensional eigenanalysis.

There are various open questions associated with this work. Although we have focused on a Gaussian sampling distribution, parts of our analysis provide sufficient conditions for general noise matrices. While qualitatively similar results should hold for sub-Gaussian distributions [5], it would be interesting to characterize how these conditions change as the tail behavior of the noise is varied away from sub-Gaussian. For instance, under bounded moment conditions, one would expect to obtain rates polynomial (as opposed to logarithmic) in the dimension $p$. It is also interesting to consider extensions of our support recovery analysis to recovery of higher rank "spiked" matrices, in the spirit of Paul and Johnstone's [32] work on $\ell_2$-approximation, as opposed to the rank-one eigenvector outer product considered here.

## APPENDIX A: MATRIX NORMS

In this appendix, we review some of the properties of matrix norms, with an emphasis on induced operator norms. Recall from (4) that for a matrix $A \in \mathbb{R}^{m \times n}$, the operator norm induced by the vector norms $\ell_p$ and $\ell_q$ (on $\mathbb{R}^m$ and $\mathbb{R}^n$, resp.) is defined by

$$\|A\|_{p,q} = \max_{\|x\|_q = 1} \|Ax\|_p \tag{56}$$

for integers $1 \leq p, q \leq \infty$. As particular examples, we have the $\ell_1$-operator norm given by $\|A\|_{1,1} = \max_{1 \leq j \leq m} \sum_{i=1}^{n} |A_{ij}|$, the $\ell_\infty$-operator norm by $\|A\|_{\infty, \infty} = \max_{1 \leq i \leq n} \sum_{j=1}^{m} |A_{ij}|$ and the spectral or $\ell_2$-operator norm by $\|A\|_{2,2} = \max\{\sigma_i(A)\}$, where $\sigma_i(A)$ are the singular values of $A$.

As a consequence of the definition (56), for any vector $x \in \mathbb{R}^n$, we have

$$\|Ax\|_p \leq \|A\|_{p,q} \|x\|_q, \tag{57}$$

a property referred to as $\|\cdot\|_{p,q}$ being *consistent* with vector norms $\|\cdot\|_p$ and $\|\cdot\|_q$ (on $\mathbb{R}^m$ and $\mathbb{R}^n$, resp.). It also follows from the definition, using (57) twice, that operator norms are consistent with themselves, in the following sense: if $A \in \mathbb{R}^{m \times n}$ and $B \in \mathbb{R}^{n \times k}$, then

$$\|AB\|_{p,q} \leq \|A\|_{p,r} \|B\|_{r,q} \tag{58}$$

for all $1 \leq p, q, r \leq \infty$.

We can also apply any vector norm to matrices, treating them as vectors, by concatenating their columns together. For example, we will use the following mixed-norm inequality

$$(59) \qquad \|AB\|_\infty \leq \|A\|_{\infty,\infty} \|B\|_\infty,$$

where $\|B\|_\infty := \max_{i,j} |B_{ij}|$ is the elementwise $\ell_\infty$-norm, and $A$ and $B$ are as defined above. For the proof, let $b_1, \ldots, b_k$ denote the columns of $B$. Then,

$$\|AB\|_\infty = \|[Ab_1, \ldots, Ab_k]\|_\infty = \max_{1 \leq i \leq k} \|Ab_i\|_\infty$$

$$\leq \|A\|_{\infty,\infty} \max_{1 \leq i \leq p} \|b_i\|_\infty = \|A\|_{\infty,\infty} \|B\|_\infty.$$

For more details, see the standard books [18, 34].

## APPENDIX B: LARGE DEVIATIONS FOR CHI-SQUARED VARIATES

The following large-deviations bounds for centralized $\chi^2$ are taken from Laurent and Massart [24]. Given a centralized $\chi^2$-variate $X$ with $d$ degrees of freedom, then for all $x \geq 0$,

$$(60a) \qquad \mathbb{P}[X - d \geq 2\sqrt{dx} + 2x] \leq \exp(-x) \quad \text{and}$$

$$(60b) \qquad \mathbb{P}[X - d \leq -2\sqrt{dx}] \leq \exp(-x).$$

We also use the following slightly different version of the bound (60a),

$$(61) \qquad \mathbb{P}[X - d \geq dx] \leq \exp(-\tfrac{3}{16}dx^2), \qquad 0 \leq x < \tfrac{1}{2},$$

due to Johnstone [19]. More generally, the analogous tail bounds for *noncentral* $\chi^2$, taken from Birgé [4], can be established via the Chernoff bound. Let $X$ be a noncentral $\chi^2$ variable with $d$ degrees of freedom and noncentrality parameter $\nu \geq 0$. Then, for all $x > 0$,

$$(62a) \qquad \mathbb{P}[X \geq (d + \nu) + 2\sqrt{(d + 2\nu)x} + 2x] \leq \exp(-x) \quad \text{and}$$

$$(62b) \qquad \mathbb{P}[X \leq (d + \nu) - 2\sqrt{(d + 2\nu)x}] \leq \exp(-x).$$

We derive here a slightly weakened but useful form of the bound (62a), valid when $\nu$ satisfies $\nu \leq Cd$ for a positive constant $C$. Under this assumption, then for any $\delta \in (0, 1)$, we have

$$(63) \qquad \mathbb{P}[X \geq (d + \nu) + 4d\sqrt{\delta}] \leq \exp\left(-\frac{\delta}{1 + 2C}d\right).$$

To establish this bound, let $x = \frac{d^2\delta}{d + 2\nu}$ for some $\delta \in (0, 1)$. From (62a), we have

$$p^* := \mathbb{P}\left[X \geq (d + \nu) + 2d\sqrt{\delta} + 2\frac{d^2}{d + 2\nu}\delta\right] \leq \exp\left(-\frac{d^2\delta}{d + 2\nu}\right)$$

$$\leq \exp\left(-\frac{\delta}{1 + 2C}d\right).$$

Moreover, we have

$$p^* \geq \mathbb{P}[X \geq (d + \nu) + 2d\sqrt{\delta} + 2d\delta] \geq \mathbb{P}[X \geq (d + \nu) + 4d\sqrt{\delta}],$$

since $\sqrt{\delta} \geq \delta$ for $\delta \in (0, 1)$.

## APPENDIX C: PROOF OF LEMMA 4

Using the form of the $\chi_n^2$ PDF, we have, for even $n$ and any $t > 0$,

$$\mathbb{P}\left[\frac{\chi_n^2}{n} > 1 + t\right]$$

$$= \frac{1}{2^{n/2}\Gamma(n/2)} \int_{(1+t)n}^{\infty} x^{n/2-1} \exp(-x/2)\, dx$$

$$= \frac{1}{2^{n/2}\Gamma(n/2)} \left\{ \frac{(n/2-1)!}{(1/2)^{(n/2-1)+1}} \exp\left(-\frac{n(1+t)}{2}\right) \sum_{i=0}^{n/2-1} \frac{1}{i!}\left(\frac{n(1+t)}{2}\right)^i \right\}$$

$$\geq \exp(-nt/2)\left[\frac{\exp(-n/2)(n/2)^{n/2-1}}{(n/2-1)!}\right](1+t)^{n/2-1},$$

where the second line uses standard integral formula (cf. Section 3.35 in the reference book [14]). Using Stirling's approximation for $(n/2-1)!$, the term within square brackets is lower bounded by $2C/\sqrt{n}$. Also, over $t \in (0, 1)$, we have $(1+t)^{-1} > 1/2$, so we conclude that

$$(64) \qquad \mathbb{P}\left[\frac{\chi_n^2}{n} > 1 + t\right] \geq \frac{C}{\sqrt{n}} \exp\left(-\frac{n}{2}[t - \log(1+t)]\right).$$

Defining the function $f(t) = \log(1+t)$, we calculate $f(0) = 0$, $f'(0) = 1$ and $f''(t) = -1/(1+t)^2$. Note that $f''(t) \geq -1$, for all $t \in \mathbb{R}$. Consequently, via a second-order Taylor series expansion, we have $f(t) - t \geq -t^2/2$. Substituting this bound into (64) yields

$$\mathbb{P}\left[\frac{\chi_n^2}{n} > 1 + t\right] \geq \frac{C}{\sqrt{n}} \exp\left(-\frac{nt^2}{2}\right)$$

as claimed.

## APPENDIX D: PROOFS FOR SECTION 4.2

**D.1. Proof of Lemma 6.** The argument we present here has a deterministic nature. In other words, we will show that if the conditions of the lemma hold for a nonrandom sequence of matrices $\Delta_{SS}$, the conclusions will follow. Thus, for example, all the references to limits may be regarded as deterministic. Then, since the conditions of the lemma are assumed to hold for a random $\Delta_{SS}$ a.a.s., it immediately follows that the conclusions hold a.a.s. To simplify the argument let us assume that $\alpha^{-1}\|\Delta_{SS}\|_{\infty,\infty} \leq \varepsilon$ for sufficiently small $\varepsilon > 0$; it turns out that $\varepsilon = \frac{1}{10}$ is enough.

We prove the lemma in steps. First, by Weyl's theorem [18, 34], eigenvalues of the perturbed matrix $\alpha z_S^* z_S^{*T} + \Delta_{SS}$ are contained in intervals of length $2\|\Delta_{SS}\|_{2,2}$ centered at eigenvalues of $\alpha z_S^* z_S^{*T}$. Since the matrix $z_S^* z_S^{*T}$ is rank one, one eigenvalue of the perturbed matrix is in the interval $[\alpha \pm \|\Delta_{SS}\|_{2,2}]$, and the remaining $k-1$ eigenvalues are in the interval $[0 \pm \|\Delta_{SS}\|_{2,2}]$. Since by assumption $2\|\Delta_{SS}\|_{2,2} \leq \alpha$ eventually, the two intervals are disjoint, and the first one contains the maximal eigenvalue $\gamma_1$ while the second contains the second largest eigenvalue $\gamma_2$. In other words, $|\gamma_1 - \alpha| \leq \|\Delta_{SS}\|_{2,2}$ and $|\gamma_2| \leq \|\Delta_{SS}\|_{2,2}$. Since $\|\Delta_{SS}\|_{2,2} \to 0$ by assumption, we conclude that $\gamma_1 \to \alpha$ and $\gamma_2 \to 0$. For the rest of the proof, take $n$ large enough so

$$|\gamma_1 \alpha^{-1} - 1| \leq \varepsilon, \tag{65}$$

where $\varepsilon > 0$ is a small number to be determined.

Now, let $\widehat{z}_S \in \mathbb{R}^k$ with $\|\widehat{z}_S\|_2 = 1$ be the eigenvector associated with $\gamma_1$, that is,

$$(\alpha z_S^* z_S^{*T} + \Delta_{SS})\widehat{z}_S = \gamma_1 \widehat{z}_S. \tag{66}$$

Taking inner products with $\widehat{z}_S$, one obtains $\alpha(z_S^{*T}\widehat{z}_S)^2 + \widehat{z}_S^T \Delta_{SS}\widehat{z}_S = \gamma_1$. Noting that $|\widehat{z}_S^T \Delta_{SS}\widehat{z}_S|$ is upper-bounded by $\|\Delta_{SS}\|_{2,2}$, we have by triangle inequality

$$|\alpha - \alpha(z_S^{*T}\widehat{z}_S)^2| = |\alpha - \gamma_1 + \gamma_1 - \alpha(z_S^{*T}\widehat{z}_S)^2|$$
$$\leq |\alpha - \gamma_1| + |\gamma_1 - \alpha(z_S^{*T}\widehat{z}_S)^2| \leq 2\|\Delta_{SS}\|_{2,2},$$

which implies $z_S^{*T}\widehat{z}_S \to 1$ (taking into account our sign convention). Take $n$ large enough so that

$$|z_S^{*T}\widehat{z}_S - 1| \leq \varepsilon \tag{67}$$

and let $u$ be the solution of

$$\alpha z_S^* + \Delta_{SS}u = \alpha u, \tag{68}$$

which is an approximation of (66) satisfied by $\widehat{z}_S$. Using triangle inequality, one has $\|u\|_\infty \leq \|z_S^*\|_\infty + \alpha^{-1}\|\Delta_{SS}\|_{\infty,\infty}\|u\|_\infty$, which implies that

$$\|u\|_\infty \leq (1 - \alpha^{-1}\|\Delta_{SS}\|_{\infty,\infty})^{-1}\|z_S^*\|_\infty \leq (1-\varepsilon)^{-1}\|z_S^*\|_\infty. \tag{69}$$

We also have

$$\|u - z_S^*\|_\infty \leq \alpha^{-1}\|\Delta_{SS}\|_{\infty,\infty}\|u\|_\infty \leq \varepsilon(1-\varepsilon)^{-1}\|z_S^*\|_\infty. \tag{70}$$

Subtracting (68) from (66), we obtain $\alpha z_S^*(z_S^{*T}\hat{z}_S-1)+\Delta_{SS}(\hat{z}_S-u)=\gamma_1\hat{z}_S-\alpha u$. Adding and subtracting $\gamma_1 u$ on the right-hand side and dividing by $\alpha$, we have

$$z_S^*(z_S^{*T}\hat{z}_S-1)+\alpha^{-1}\Delta_{SS}(\hat{z}_S-u)=\gamma_1\alpha^{-1}(\hat{z}_S-u)+(\gamma_1\alpha^{-1}-1)u,$$

which implies

$$\|\hat{z}_S-u\|_\infty \le (|\gamma_1\alpha^{-1}|-\alpha^{-1}\|\!|\Delta_{SS}|\!\|_{\infty,\infty})^{-1}$$
$$\times\{|z_S^{*T}\hat{z}_S-1|\cdot\|z_S^*\|_\infty+|\gamma_1\alpha^{-1}-1|\cdot\|u\|_\infty\}$$
$$\le (1-2\varepsilon)^{-1}[\varepsilon+\varepsilon(1-\varepsilon)^{-1}]\cdot\|z_S^*\|_\infty,$$

where the last inequality follows from (65), (67) and (69). Combining with the bound (70) on $\|u-z_S^*\|_\infty$ yields

$$\frac{\|\hat{z}_S-z_S^*\|_\infty}{\|z_S^*\|_\infty}\le\frac{\varepsilon}{1-2\varepsilon}+\frac{\varepsilon}{(1-2\varepsilon)(1-\varepsilon)}+\frac{\varepsilon}{1-\varepsilon}$$
$$\le\frac{3\varepsilon}{(1-2\varepsilon)^2}.$$

Finally, we take $\varepsilon=\frac{1}{10}$ to conclude $\|\hat{z}_S-z_S^*\|_\infty\le\frac{1}{2}\|z_S^*\|_\infty=\frac{1}{2\sqrt{k}}$ a.a.s., as claimed.

**D.2. Proof of Lemma 7.** Recall that by definition, $\tilde{z}_S=\hat{z}_S/\|\hat{z}_S\|_1$. Using the identity $\text{sign}(\hat{z}_S)^T\hat{z}_S=\|\hat{z}_S\|_1$ yields $\hat{U}_{S^cS}\hat{z}_S=\rho_n^{-1}\Delta_{S^cS}\tilde{z}_S$, which is the desired equation. It only remains to prove that $\hat{U}_{S^cS}$ is indeed a valid sign matrix.

First note that from (28) we have $|\hat{z}_i|\in[\frac{1}{2\sqrt{k}},\frac{3}{2\sqrt{k}}]$ for $i\in S$, which implies that $\|\hat{z}_S\|_1\in[\frac{\sqrt{k}}{2},\frac{3\sqrt{k}}{2}]$. Thus, $\|\tilde{z}_S\|_2=1/(\|\hat{z}_S\|_1)\le\frac{2}{\sqrt{k}}$. Now we can write

$$\max_{i\in S^c,j\in S}|\hat{U}_{ij}|\le\rho_n^{-1}\|\Delta_{S^cS}\tilde{z}_S\|_\infty$$
$$\le\rho_n^{-1}\|\!|\Delta_{S^cS}|\!\|_{\infty,2}\|\tilde{z}_S\|_2$$
$$\le\frac{2k}{\beta}\frac{\delta}{\sqrt{k}}\frac{2}{\sqrt{k}}$$
$$=\frac{4}{\beta}\delta,$$

so that taking $\delta\le\frac{\beta}{4}$ completes the proof.

**D.3. Proof of Lemma 8.** Here we provide the proof for the case $\Gamma_{p-k}=I_{p-k}$; necessary modifications for the general case are discussed in Section

**4.4.** First, let us bound the cross-term in (32). Recall that $\widetilde{z}_S = \widehat{z}_S / \|\widehat{z}_S\|_1$. Also, by our choice (30) of $\widehat{U}_{S^c S}$, we have

$$\Phi_{S^c S} = \Delta_{S^c S} - \rho_n \widehat{U}_{S^c S} = \Delta_{S^c S} - \Delta_{S^c S} \widetilde{z}_S \operatorname{sign}(\widehat{z}_S)^T.$$

Now, using sub-multiplicative property of operator norms [see relation (58) in Appendix A], we can write

$$
\begin{aligned}
(71) \qquad \|\Phi_{S^c S}\|_{\infty,2} &= \|\Delta_{S^c S}(I_{p-k} - \widetilde{z}_S \operatorname{sign}(\widehat{z}_S)^T)\|_{\infty,2} \\
&\leq \|\Delta_{S^c S}\|_{\infty,2} \cdot \|(I_{p-k} - \widetilde{z}_S \operatorname{sign}(\widehat{z}_S)^T)\|_{2,2} \\
&\leq \|\Delta_{S^c S}\|_{\infty,2} \cdot (1 + |\widetilde{z}_{S_2} \operatorname{sign}(\widehat{z}_S)^T \widetilde{z}_S|) \\
&\leq 3\|\Delta_{S^c S}\|_{\infty,2},
\end{aligned}
$$

where we have also used the fact that $\|ab^T\|_{2,2} = \|a\|_2 \|b\|_2$, and $\|\widetilde{z}_S\|_2 = 1/(\|\widehat{z}_S\|_1) \leq \frac{2}{\sqrt{k}}$, using the bound (28). Recall the decomposition $x = (u, v)$, where $u = \mu \widehat{z}_S + \widehat{z}_S^\perp$ with $\mu^2 + \|\widehat{z}_S^\perp\|_2^2 \leq 1$. Also, by our choice (30) of $\widehat{U}_{S^c S}$, we have $\Phi_{S^c S} u = \Phi_{S^c S} \widehat{z}_S^\perp$. Thus,

$$
\begin{aligned}
(72) \qquad \max_u |2 v^T \Phi_{S^c S} u| &\leq \max_{\substack{\|\tilde{u}\|_2 \leq \sqrt{1-\mu^2}, \\ \tilde{u} \perp z_S}} |2 v^T \Phi_{S^c S} \tilde{u}| \\
&\leq \sqrt{1-\mu^2} \max_{\|\tilde{u}\|_2 \leq 1} |2 v^T \Phi_{S^c S} \tilde{u}|.
\end{aligned}
$$

Using Hölder's inequality, we have

$$
\begin{aligned}
(73) \qquad \max_{\|\tilde{u}\|_2 \leq 1} |2 v^T \Phi_{S^c S} \tilde{u}| &\leq 2\|v\|_1 \max_{\|\tilde{u}\|_2 \leq 1} \|\Phi_{S^c S} \tilde{u}\|_\infty \\
&\leq 2\|v\|_1 \|\Phi_{S^c S}\|_{\infty,2} \\
&\leq 6\|v\|_1 \frac{\delta}{\sqrt{k}},
\end{aligned}
$$

where we have used bound (71) and applied condition (31). We now turn to the last term in the decomposition (32), namely $v^T \Phi_{S^c S^c} v = v^T \Delta_{S^c S^c} v - \rho_n v^T \widehat{U}_{S^c S^c} v$. In order to minimize this term, we use our freedom to choose $\widehat{U}_{S^c S^c}(x) = \operatorname{sign}(v) \operatorname{sign}(v)^T$, so that $-\rho_n v^T \widehat{U}_{S^c S^c} v$ simply becomes $-\rho_n \|v\|_1^2$.

Define the objective function $f^* := \max_x x^T \Phi x$. Also let $H = n^{-1/2} G_{S^c}$, where $G_{S^c} = (G_{ij})$ for $1 \leq i \leq n$ and $j \in S^c$. Noting that $\Delta_{S^c S^c} = H^T H - I_m$ (with $m = p - k$) and using the bounds (33), (72) and (73), we obtain the following bound on the objective:

$$
\begin{aligned}
f^* &\leq \max_u u^T \Phi_{SS} u + \max_{u,v} 2 v^T \Phi_{S^c S} u + \max_v v^T \Phi_{S^c S^c} v \\
(74) \qquad &\leq [\mu^2 \gamma_1 + (1-\mu^2)\gamma_2]
\end{aligned}
$$

$$+ (1-\mu^2) \underbrace{\left[ \max_{\|v\|_2 \le 1} \left\{ 6\|v\|_1 \frac{\delta}{\sqrt{k}} + \|Hv\|_2^2 - \|v\|_2^2 - \rho_n \|v\|_1^2 \right\} \right]}_{g^*}.$$

In obtaining the last inequality, we have used the change of variable $v \to (\sqrt{1-\mu^2})v$, with some abuse of notation, and exploited the inequality $\|v\|_2 \le \sqrt{1-\mu^2}$. (Note that this bound follows from the identity $\|x\|_2^2 = 1 = \mu^2 + \|\hat{z}_S^\perp\|_2^2 + \|v\|_2^2$.)

Let $v^*$ be the optimal solution to problem $g^*$ in (74); note that it is random due to the presence of $H$. Also, set $\mathbb{S} = \{(\eta_{ij}, \ell_{ij})\}$ where $i$ and $j$ range over $\{1, 2, \dots, \lceil \sqrt{m} \rceil\}$ and

$$\eta_{ij} = \frac{i}{\sqrt{m}}, \qquad \ell_{ij} = \frac{i}{\sqrt{m}} j.$$

Note that $\mathbb{S}$ satisfies the condition of the lemma, namely $|\mathbb{S}| = \lceil \sqrt{m} \rceil^2 = \mathcal{O}(m)$.

Since $\|v^*\|_2 \le 1$, and $\|v^*\|_2 \le \|v^*\|_1 \le \sqrt{m}\|v^*\|_2$, there exists[3] $(\eta^*, \ell^*) \in \mathbb{S}$ such that

$$\eta^* - \frac{1}{\sqrt{m}} < \|v^*\|_2 \le \eta^*,$$

$$\ell^* - 3 < \|v^*\|_1 \le \ell^*.$$

Thus, using condition (34), we have

$$\|Hv^*\|_2 \le \max_{\substack{\|v\|_2 \le \eta^*, \\ \|v\|_1 \le \ell^*}} \|Hv\|_2 \le \eta^* + \frac{\delta}{\sqrt{k}} \ell^* + \varepsilon$$

$$\le \|v^*\|_2 + \frac{1}{\sqrt{m}} + \frac{\delta}{\sqrt{k}}(\|v^*\|_1 + 3) + \varepsilon.$$

To simplify notation, let

$$(75) \qquad A = A(\varepsilon, \delta, m, k) := 1/\sqrt{m} + 3\delta/\sqrt{k} + \varepsilon,$$

so that the bound in the above display may be written as $\|v^*\|_2 + \delta \|v^*\|_1/\sqrt{k} + A$. Now, we have

$$\|Hv^*\|_2^2 - \|v^*\|_2^2 \le 2\|v^*\|_2 \left( \delta \frac{\|v^*\|_1}{\sqrt{k}} + A \right) + \left( \delta \frac{\|v^*\|_1}{\sqrt{k}} + A \right)^2$$

---

[3]Let $i^* = \lceil \sqrt{m} \|v^*\|_2 \rceil$ and $\eta^* = \frac{i^*}{\sqrt{m}}$. Using the fact that, for any $x \in \mathbb{R}$, $\lceil x \rceil - 1 < x \le \lceil x \rceil$, we have $\eta^* - 1/\sqrt{m} < \|v^*\|_2 \le \eta^*$ or, equivalently, $\|v^*\|_2 = \eta^* + \xi$ where $-1/\sqrt{m} < \xi \le 0$. Now let $j^* = \lceil \frac{\|v^*\|_1}{\|v^*\|_2} \rceil$. One has $(j^* - 1)\|v^*\|_2 < \|v^*\|_1 \le j^* \|v^*\|_2$ which, using the fact that $\|v^*\|_2 \le 1$, implies $j^* \|v^*\|_2 - 1 < \|v^*\|_1 \le j^* \|v^*\|_2$. This in turn implies

$$j^* \eta^* + j^* \xi - 1 < \|v^*\|_1 \le j^* \eta^*.$$

Take $\ell^* = j^* \eta^*$ and note that $j^* \xi - 1 > -3$, since $j^*$ is at most $\lceil \sqrt{m} \rceil$.

(76)
$$\leq 2\left(\delta \frac{\|v^*\|_1}{\sqrt{k}} + A\right) + \left(\delta \frac{\|v^*\|_1}{\sqrt{k}} + A\right)^2.$$

Using this in (74) and recalling from (26) that $\rho_n = \beta/(2k)$, we obtain the following bound:

$$g^* \leq 6\delta \frac{\|v^*\|_1}{\sqrt{k}} + 2\left(\delta \frac{\|v^*\|_1}{\sqrt{k}} + A\right) + \left(\delta \frac{\|v^*\|_1}{\sqrt{k}} + A\right)^2 - \frac{\beta}{2}\left(\frac{\|v^*\|_1}{\sqrt{k}}\right)^2.$$

Note that this is quadratic in $\|v^*\|_1/\sqrt{k}$, that is,

$$g^* \leq a\left(\frac{\|v^*\|_1}{\sqrt{k}}\right)^2 + b\left(\frac{\|v^*\|_1}{\sqrt{k}}\right) + c,$$

where

$$a = \delta^2 - \frac{\beta}{2}, \qquad b = 8\delta + 2\delta A \quad \text{and} \quad c = 2A + A^2.$$

By choosing $\delta$ sufficiently small, say $\delta^2 \leq \beta/4$, we can make $a$ negative. This makes the quadratic form $ax^2 + bx + c$ achieve a maximum of $c + b^2/4(-a)$, at the point $x^* = b/2(-a)$. Note that we have $b/2(-a) \to 0$ and $c \to 0$ as $\varepsilon, \delta \to 0$ and $m, k \to \infty$. Consequently, we can make this maximum (and hence $g^*$) arbitrarily small eventually, say less than $\alpha/2$, by choosing $\delta$ and $\varepsilon$ sufficiently small.

Combining this bound on $g^*$ with our bound (74) on $f^*$, and recalling that $\gamma_1 \to \alpha$ and $\gamma_2 \to 0$ by Lemma 6, we conclude that

$$f^* \leq \mu^2(\alpha + o(1)) + (1 - \mu^2)\left[\frac{\alpha}{2} + o(1)\right] \leq \alpha + o(1)$$

as claimed.

## APPENDIX E: PROOF OF LEMMA 9

In this appendix, we prove Lemma 9, a general result on $\|\cdot\|_{\infty,2}$-norm of Wishart matrices. Some of the intermediate results are of independent interest and are stated as separate lemmas. Two sets of large deviation inequalities will be used, one for chi-squared RVs $\chi_n^2$ and one for "sums of Gaussian product" random variates. To define the latter precisely, let $Z_1$ and $Z_2$ be independent Gaussian RVs, and consider the sum $\sum_{i=1}^n X_i$ where $X_i \overset{\text{i.i.d.}}{\sim} Z_1 Z_2$, for $1 \leq i \leq n$. The following tail bounds are known [4, 21]:

(77)
$$\mathbb{P}\left(\left|n^{-1} \sum_{i=1}^n X_i\right| > t\right) \leq C \exp(-3nt^2/2) \qquad \text{as } t \to 0;$$

(78)
$$\mathbb{P}(|n^{-1}\chi_n^2 - 1| > t) \leq 2\exp(-3nt^2/16), \qquad 0 \leq t < 1/2,$$

where $C$ is some positive constant.

Let $W$ be a $p \times p$ centered Wishart matrix as defined in (37). Consider the following linear combination of off-diagonal entries of the first row:

$$\sum_{j=2}^{n} a_j W_{1j} = n^{-1} \sum_{i=1}^{n} g_1^i \sum_{j=2}^{p} g_j^i a_j.$$

Let $\xi^i := \|a\|_2^{-1} \sum_{j=2}^{p} g_j^i a_j$, where $a = (a_2, \ldots, a_p) \in \mathbb{R}^{p-1}$. Note that $\{\xi^i\}_{i=1}^{n}$ is a collection of independent standard Gaussian RVs. Moreover, $\{\xi^i\}_{i=1}^{n}$ is independent of $\{g_1^i\}_{i=1}^{n}$. Now we have

$$\sum_{j=2}^{p} a_j W_{1j} = n^{-1} \|a\|_2 \sum_{i=1}^{n} g_1^i \xi^i,$$

which is a (scaled) sum of Gaussian products (as defined above). Using (77), we obtain

$$(79) \qquad \mathbb{P}\left( \left| \sum_{j=2}^{p} a_j W_{1j} \right| > t \right) \leq C \exp(-3nt^2/2\|a\|_2^2).$$

Combining the bounds in (79) and (78), we can bound a full linear combination of first-row entries. More specifically, let $x = (x_1, \ldots, x_p) \in \mathbb{R}^p$, with $x_1 \neq 0$ and $\sum_{j=2}^{p} x_j \neq 0$, and consider the linear combination $\sum_{j=1}^{p} x_j W_{1j}$. Noting that $W_{11} = n^{-1} \sum_i (g_1^i)^2 - 1$ is a centered $\chi_n^2$, we obtain

$$\mathbb{P}\left[ \left| \sum_{j=1}^{p} x_j W_{1j} \right| > t \right] \leq \mathbb{P}\left( |x_1 W_{11}| + \left| \sum_{j=2}^{p} x_j W_{1j} \right| > t \right)$$

$$\leq \mathbb{P}[|x_1 W_{11}| > t/2] + \mathbb{P}\left[ \left| \sum_{j=2}^{p} x_j W_{1j} \right| > t/2 \right]$$

$$\leq 2 \exp\left( -\frac{3nt^2}{16 \cdot 4x_1^2} \right) + C \exp\left( -\frac{3nt^2}{2 \cdot 4 \sum_{j=2}^{p} x_j^2} \right)$$

$$\leq 2 \max\{2, C\} \exp\left( -\frac{3nt^2}{16 \cdot 4 \sum_{j=1}^{p} x_j^2} \right).$$

Note that the last inequality holds, in general, for $x \neq 0$. Since there is nothing special about the "first" row, we can conclude the following.

LEMMA 15. *For $t > 0$ small enough, there are (numerical constants) $c > 0$ and $C > 0$ such that for all $x \in \mathbb{R}^p \setminus \{0\}$,*

$$(80) \qquad \mathbb{P}\left( \left| \sum_{j=1}^{p} x_j W_{ij} \right| > t \right) \leq C \exp(-cnt^2/\|x\|_2^2)$$

*for $1 \leq i \leq p$.*

Now, let $\mathcal{I}, \mathcal{J} \subset \{1, \ldots, p\}$ be index sets,[4] both allowed to depend on $p$ (though we have omitted the dependence for brevity). Choose $x$ such that $x_j = 0$ for $j \notin \mathcal{J}$ and $\|x_{\mathcal{J}}\|_2 = 1$. Note that $\|W_{\mathcal{I}, \mathcal{J}} x_{\mathcal{J}}\|_\infty = \max_{i \in \mathcal{I}} |\sum_{j \in \mathcal{J}} W_{ij} x_j| = \max_{i \in \mathcal{I}} |\sum_{j=1}^p W_{ij} x_j|$, suggesting the following lemma.

LEMMA 16. *Consider some index set $\mathcal{I}$ such that $|\mathcal{I}| \to \infty$ and $n^{-1} \log |\mathcal{I}| \to 0$ as $n, p \to \infty$, and some $x_{\mathcal{J}} \in S^{|\mathcal{J}|-1}$. Then, there exists an absolute constant $B > 0$ such that*

$$\|W_{\mathcal{I}, \mathcal{J}} x_{\mathcal{J}}\|_\infty \le B \sqrt{\frac{\log |\mathcal{I}|}{n}} \tag{81}$$

*as $n, p \to \infty$, with probability 1.*

PROOF. Applying the union bound in conjunction with the bound (80) yields

$$\mathbb{P}\left(\max_{i \in \mathcal{I}} \left|\sum_{j \in \mathcal{J}} W_{ij} x_j\right| > t\right) \le |\mathcal{I}| C \exp(-cnt^2). \tag{82}$$

Letting $t = B\sqrt{n^{-1} \log |\mathcal{I}|}$, the right-hand side simplifies to $C \exp(-(cB^2 - 1) \log |\mathcal{I}|)$. Taking $B > \sqrt{2c^{-1}}$ and applying Borel–Cantelli lemma completes the proof. □

Note that as a corollary, setting $x_{\mathcal{J}} = (1, 0, \ldots, 0)$ yields bounds on the $\infty$-norm of columns (or, equivalently, rows) of Wishart matrices.

Lemma 16 may be used to obtain the desired bound on $\|W_{\mathcal{I}, \mathcal{J}}\|_{\infty, 2}$. For simplicity, let $y \in \mathbb{R}^{|\mathcal{J}|}$ represent a generic $|\mathcal{J}|$-vector. Recall that $\|W_{\mathcal{I}, \mathcal{J}}\|_{\infty, 2} = \max_{y \in S^{|\mathcal{J}|-1}} \|W_{\mathcal{I}, \mathcal{J}} y\|_\infty$. We use a standard discretization argument, covering the unit $\ell^2$-ball of $\mathbb{R}^{|\mathcal{J}|}$ using an $\varepsilon$-net, say $\mathcal{N}$. It can be shown [27] that there exists such a net with cardinality $|\mathcal{N}| < (3/\varepsilon)^{|\mathcal{J}|}$. For every $y \in S^{|\mathcal{J}|-1}$, let $u_y \in \mathcal{N}$ be the point such that $\|y - u_y\|_2 \le \varepsilon$. Then

$$\|W_{\mathcal{I}, \mathcal{J}} y\|_\infty \le \|W_{\mathcal{I}, \mathcal{J}}\|_{\infty, 2} \|y - u_y\|_2 + \|W_{\mathcal{I}, \mathcal{J}} u_y\|_\infty$$
$$\le \|W_{\mathcal{I}, \mathcal{J}}\|_{\infty, 2} \varepsilon + \|W_{\mathcal{I}, \mathcal{J}} u_y\|_\infty.$$

Taking the maximum over $y \in S^{|\mathcal{J}|-1}$ and rearranging yields the inequality

$$\|W_{\mathcal{I}, \mathcal{J}}\|_{\infty, 2} \le (1 - \varepsilon)^{-1} \max_{u \in \mathcal{N}} \|W_{\mathcal{I}, \mathcal{J}} u\|_\infty. \tag{83}$$

---

[4] We always assume that these index sets form an increasing sequence of sets. More precisely, with $\mathcal{I} = \mathcal{I}_p$, we assume $\mathcal{I}_1 \subset \mathcal{I}_2 \subset \cdots$. We also assume $|\mathcal{I}_p| \to \infty$ as $p \to \infty$.

Using this bound (83), we can now provide the proof of Lemma 9 as follows. Let $\mathcal{N} = \{u_1, \dots, u_{|\mathcal{N}|}\}$ be a $\frac{1}{2}$-net of the ball $S^{|\mathcal{J}|-1}$, with cardinality $|\mathcal{N}| < 6^{|\mathcal{J}|}$. Then, from our bound (83), we have

$$\mathbb{P}(\|W_{\mathcal{I},\mathcal{J}}\|_{\infty,2} > t) \leq \mathbb{P}\Big(2\max_{u \in \mathcal{N}} \|W_{\mathcal{I},\mathcal{J}}u\|_{\infty} > t\Big)$$
$$\leq |\mathcal{N}| \cdot \mathbb{P}(\|W_{\mathcal{I},\mathcal{J}}u_1\|_{\infty} > t/2)$$
$$\leq 6^{|\mathcal{J}|} \cdot C|\mathcal{I}| \exp(-cnt^2/4).$$

In the last line, we used (82). Taking $t = D'' \frac{\sqrt{|\mathcal{J}|} + \sqrt{\log |\mathcal{I}|}}{\sqrt{n}}$ with $D''$ large enough and using Borel–Cantelli lemma completes the proof.

## APPENDIX F: PROOF OF LEMMA 14

The mixture covariance can be expressed as

$$\Sigma_M := \mathbb{E}[xx^T] = \mathbb{E}[\mathbb{E}[xx^T|U]]$$
$$= \sum_{S \in \widetilde{\mathbb{S}}} \frac{1}{|\widetilde{\mathbb{S}}|} \mathbb{E}[xx^T|U = S]$$
$$= \sum_{S \in \widetilde{\mathbb{S}}} \frac{1}{|\widetilde{\mathbb{S}}|} (I_p + \beta z^*(S)z^*(S)^T)$$
$$= I_p + \frac{\beta}{|\widetilde{\mathbb{S}}|} \sum_{S \in \widetilde{\mathbb{S}}} z^*(S)z^*(S)^T =: I_p + \frac{\beta}{k|\widetilde{\mathbb{S}}|} Y,$$

where

$$Y_{ij} = \sum_{S \in \widetilde{\mathbb{S}}} [\sqrt{k}z^*(S)]_i [\sqrt{k}z^*(S)^T]_j = \sum_{S \in \widetilde{\mathbb{S}}} \mathbb{1}\{i \in S\}\mathbb{1}\{j \in S\}$$
$$= \sum_{S \in \widetilde{\mathbb{S}}} \mathbb{1}\{\{i,j\} \subset S\}.$$

Let $R := \{1, \dots, k-1\}$ and $R^c := \{k, \dots, p\}$. Note that we always have $R \subset S$ for $S \in \widetilde{\mathbb{S}}$. In general, we have

$$Y_{ij} = \begin{cases} |\widetilde{\mathbb{S}}|, & \text{if both } i, j \in R, \\ 1, & \text{if exactly one of } i \text{ or } j \in R, \\ 0, & \text{if both } i, j \notin R. \end{cases}$$

Consequently, $Y$ takes the form

$$Y = \begin{pmatrix} |\widetilde{\mathbb{S}}| & \cdots & |\widetilde{\mathbb{S}}| & 1 & 1 & \cdots & 1 \\ \vdots & \ddots & \vdots & \vdots & \vdots & \ddots & \vdots \\ |\widetilde{\mathbb{S}}| & \cdots & |\widetilde{\mathbb{S}}| & 1 & 1 & \cdots & 1 \\ \hline 1 & \cdots & 1 & 1 & 0 & \cdots & 0 \\ 1 & \cdots & 1 & 0 & 1 & \cdots & 0 \\ \vdots & \ddots & \vdots & \vdots & \vdots & \ddots & \vdots \\ 1 & \cdots & 1 & 0 & 0 & \cdots & 1 \end{pmatrix} \quad \text{or} \quad Y = \begin{pmatrix} |\widetilde{\mathbb{S}}|\vec{1}_R\vec{1}_R^T & \vec{1}_R\vec{1}_{R^c}^T \\ \vec{1}_{R^c}\vec{1}_R^T & I_{R^c \times R^c} \end{pmatrix},$$

where $\vec{1}_R$, for example, denotes the vector of all ones over the index set $R$. We conjecture an eigenvector of the form

$$v = \begin{pmatrix} \vec{1}_R \\ b\vec{1}_{R^c} \end{pmatrix}$$

and let us denote the associated eigenvalue as $\lambda$. Thus, we assume $Yv = \lambda v$, or, in more detail,

$$|\widetilde{\mathbb{S}}||R|\vec{1}_R + b|R^c|\vec{1}_R = \lambda\vec{1}_R,$$

$$|R|\vec{1}_{R^c} + b\vec{1}_{R^c} = \lambda b\vec{1}_{R^c},$$

where we have used, for example, $\vec{1}_R^T\vec{1}_R = |R|$. Note that $|R^c| = |\widetilde{\mathbb{S}}| = p - k + 1$. Rewriting in terms of $|\widetilde{\mathbb{S}}|$, we get

$$|\widetilde{\mathbb{S}}|(|R| + b) = \lambda,$$

$$|R| + b = \lambda b$$

from which we conclude, assuming $\lambda \neq 0$, that $b = \frac{1}{|\widetilde{\mathbb{S}}|}$. This, in turn, implies $\lambda = |\widetilde{\mathbb{S}}||R| + 1$.

Thus far, we have determined an eigenpair. We can now subtract $\lambda(v/\|v\|_2)(v/\|v\|_2)^T = (\lambda/\|v\|_2^2)vv^T$ and search for the rest of the eigenvalues in the remainder. Note that

$$\frac{\lambda}{\|v\|_2^2} = \frac{\lambda}{|R| + b^2|R^c|} = \frac{|\widetilde{\mathbb{S}}||R| + 1}{|R| + |\widetilde{\mathbb{S}}|^{-1}} = |\widetilde{\mathbb{S}}|.$$

Thus, we have

$$\frac{\lambda}{\|v\|_2^2}vv^T = \begin{pmatrix} |\widetilde{\mathbb{S}}|\vec{1}_R\vec{1}_R^T & \vec{1}_R\vec{1}_{R^c}^T \\ \vec{1}_{R^c}^T\vec{1}_R & \frac{1}{|\widetilde{\mathbb{S}}|}\vec{1}_{R^c}\vec{1}_{R^c}^T \end{pmatrix}$$

implying

$$Y - \frac{\lambda}{\|v\|_2^2}vv^T = \begin{pmatrix} 0 & 0 \\ 0 & I - \frac{1}{|\widetilde{\mathbb{S}}|}\vec{1}_{R^c}\vec{1}_{R^c}^T \end{pmatrix}.$$

The nonzero block of the remainder has one eigenvalue equal to $1 - \frac{|R^c|}{|\mathbb{S}|} = 0$ and the rest of $|R^c| - 1$ of its eigenvalues equal to 1. Thus, the remainder has $|R| + 1$ of its eigenvalues equal to zero and $|R^c| - 1$ of them equal to one.

Overall, we conclude that eigenvalues of $Y$ are as follows:

$$\begin{cases} |\tilde{\mathbb{S}}||R| + 1, & 1 \text{ time,} \\ 1, & |R^c| - 1 \text{ times,} \\ 0, & |R| \text{ times} \end{cases}$$

or

$$\begin{cases} (p - k + 1)(k - 1) + 1, & 1 \text{ time,} \\ 1, & p - k \text{ times,} \\ 0, & k - 1 \text{ times.} \end{cases}$$

The eigenvalues of $Y$ are mapped to those of $\Sigma_M$ by the affine map $x \to 1 + \frac{\beta}{k|\mathbb{S}|} x$, so that $\Sigma_M$ has eigenvalues

$$(84) \qquad 1 + \frac{\beta(k-1)}{k} + \frac{\beta}{k(p-k+1)}, \qquad 1 + \frac{\beta}{k(p-k+1)}, \qquad 1$$

with multiplicities 1, $p - k$ and $k - 1$, respectively. The log determinant stated in the lemma then follows by straightforward calculation.

## APPENDIX G: PROOF OF THEOREM 2(A)

Since in part (a) of the theorem we are using the weaker scaling $n > \theta_{\mathrm{wr}} k^2 \log(p - k)$, we have more freedom in choosing the sign matrix $\widehat{U}$. We choose the upper-left block $\widehat{U}_{SS}$ as in part (b) so that Lemma 6 applies. Also let $\widehat{z} := (\widehat{z}_S, \vec{0}_{S^c})$ as in (29), where $\widehat{z}_S$ is the (unique) maximal eigenvector of the $k \times k$ block $\Phi_{SS}$; it has the correct sign by Lemma 6. We set the off-diagonal and lower-right blocks of the sign matrix to

$$(85) \qquad \widehat{U}_{S^c S} = \frac{1}{\rho_n} \Delta_{S^c S}, \qquad \widehat{U}_{S^c S^c} = \frac{1}{\rho_n} \Delta_{S^c S^c},$$

so that $\Phi_{S^c S} = 0$ and $\Phi_{S^c S^c} = 0$. With these blocks of $\Phi$ being zero, $\widehat{z}$ is the maximal eigenvector of $\Phi$, hence an optimal solution of (13), if and only if $\widehat{z}_S$ is the maximal eigenvector of $\Phi_{SS}$; the latter is true by definition. Note that this argument is based on the remark following Lemma 5. It only remains to show that the choices of (85) lead to valid sign matrices.

Recalling that vector $\infty$-norm of a matrix $A$ is $\|A\|_\infty := \max_{i,j} |A_{i,j}|$ (see Appendix A), we need to show $\|\widehat{U}_{S^c S}\|_\infty \leq 1$ and $\|\widehat{U}_{S^c S^c}\|_\infty \leq 1$. Using the notation of Section 4.4 and the mixed-norm inequality (59), we have

$$\|\widehat{U}_{S^c S}\|_\infty = \frac{\sqrt{\beta}}{\rho_n} \|\tilde{h}_{S^c} z_S^{*T}\|_\infty \leq \frac{\sqrt{\beta}}{\rho_n} \|\tilde{h}_{S^c}\|_{\infty,\infty} \|z_S^{*T}\|_\infty$$

$$= \frac{\sqrt{\beta}}{\rho_n} \|\tilde{h}_{S^c}\|_\infty \|z_S^*\|_\infty$$

$$\leq \frac{\sqrt{\beta}}{\rho_n} \|\|\Gamma_{p-k}^{1/2}\|\|_{\infty,\infty} \|h_{S^c}\|_\infty \|z_S^*\|_\infty$$

$$= \frac{2k}{\sqrt{\beta}} \mathcal{O}(1) \mathcal{O}\left(\sqrt{\frac{\log(p-k)}{n}}\right) \frac{1}{\sqrt{k}} = \mathcal{O}(1) \frac{1}{\sqrt{k}} \to 0,$$

where the last line follows under the scaling assumed and assumption (6a) on $\|\|\Gamma_{p-k}^{1/2}\|\|_{\infty,\infty}$. For the lower-right block, we use the mixed-norm inequality (59) twice together with symmetry to obtain

$$\|\widehat{U}_{S^c S^c}\|_\infty = \frac{1}{\rho_n} \|\widetilde{W}_{S^c S^c}\|_\infty = \frac{1}{\rho_n} \|\Gamma_{p-k}^{1/2} W_{S^c S^c} \Gamma_{p-k}^{1/2}\|_\infty$$

$$\leq \frac{1}{\rho_n} \|\|\Gamma_{p-k}^{1/2}\|\|_{\infty,\infty}^2 \|W_{S^c S^c}\|_\infty$$

$$= \frac{2k}{\beta} \mathcal{O}(1) \mathcal{O}\left(\sqrt{\frac{\log(p-k)}{n}}\right),$$

which can be made less than one by choosing $\theta_{\mathrm{wr}}$ large enough. The bound on $\|W_{S^c S^c}\|_\infty$ used in the last line can be obtained using arguments similar to those of Lemma 9. The proof is complete.

**Acknowledgments.** We thank Alexandre d'Aspremont and Laurent El Ghaoui for helpful discussions, and the anonymous reviewers for their careful reading and helpful suggestions.

DEPARTMENT OF ELECTRICAL ENGINEERING
AND COMPUTER SCIENCE
UNIVERSITY OF CALIFORNIA, BERKELEY
BERKELEY, CALIFORNIA 94720
USA
E-MAIL: amini@eecs.berkeley.edu

DEPARTMENT OF STATISTICS
UNIVERSITY OF CALIFORNIA, BERKELEY
BERKELEY, CALIFORNIA 94720
USA
E-MAIL: wainwrig@stat.berkeley.edu